%% file: main.tex
\documentclass{article}
\usepackage[a4paper, total={16.5cm, 26cm}]{geometry}
\usepackage[table]{xcolor}

\usepackage[scr=rsfs]{mathalpha}

\usepackage{tikz}
\usetikzlibrary{}
\usetikzlibrary{matrix,chains,positioning,decorations.pathreplacing,arrows}
\usetikzlibrary{positioning,calc}
\usetikzlibrary{graphs,graphs.standard,quotes}

\usepackage{pgfplots}
\pgfplotsset{compat=1.17} 
\usepackage[utf8]{inputenc}
\usepackage[square,sort,comma,numbers]{natbib}
\usepackage{graphicx}
\graphicspath{ {./images/} }
\usepackage{amsmath}
\usepackage{amsfonts}
\usepackage{amssymb}
\usepackage{algorithm}
\usepackage{algpseudocode}
\usepackage{mathtools}
\usepackage{mathrsfs}
\usepackage{tabularx}
\usepackage{multirow}
\usepackage{float}
\usepackage{subcaption}
\usepackage[font={footnotesize}, margin=0.5cm]{caption}
\usepackage{stackengine}
\usepackage{url}
\usepackage{booktabs}
\usepackage{makecell}
\usepackage{lineno}
\usepackage{bbm}
\usepackage[bookmarks,hidelinks]{hyperref}
\usepackage{accents}

\newcommand{\norm}[1]{\left\|#1\right\|}
\newcommand{\normdg}[1]{{\left\vert\kern-0.25ex\left\vert\kern-0.25ex\left\vert #1 
    \right\vert\kern-0.25ex\right\vert\kern-0.25ex\right\vert}}
\newcommand{\tdev}[2]{\frac{\mathrm{d}#1}{\mathrm{d}#2}}

\newcommand{\bref}[1]{(\ref{#1})}

\newcommand{\mbf}[1]{\boldsymbol{#1}}

\algrenewcommand\algorithmiccomment[2][\footnotesize]{{#1\hfill\(\triangleright\) \textit{#2}}}

\newlength{\bibitemsep}
\newlength{\bibparskip}
\let\oldthebibliography\thebibliography
\renewcommand\thebibliography[1]{%
  \footnotesize
  \oldthebibliography{#1}%
  \setlength{\parskip}{\bibitemsep}%
  \setlength{\itemsep}{\bibparskip}%
}

\title{
    Neural ordinary differential equations for model order reduction of stiff systems
}
\date{\today\vspace{-5ex}}
\author{Matteo Caldana$^{a,}$\thanks{Corresponding author: {\tt matteo.caldana@polimi.it}},
Jan S. Hesthaven$^b$\\[0.3cm]
\small\textit{$^a$MOX, Dipartimento di Matematica, 
Politecnico di Milano,
}\\ \small\textit{
Piazza Leonardo da Vinci 32,
20133 Milano, Italy
}\\
\small\textit{$^b$Karlsruhe Institute of Technology,}\\ \small\textit{Kaiserstrasse 12,
76131 Karlsruhe, Germany}}

\begin{document}

\maketitle
{\subsection*{\centering Abstract}
\small 
\begin{changemargin}{1cm}{1cm}
    Neural Ordinary Differential Equations (ODEs) represent a significant advancement at the intersection of machine learning and dynamical systems, offering a continuous-time analog to discrete neural networks. Despite their promise, deploying neural ODEs in practical applications often encounters the challenge of stiffness, a condition where rapid variations in some components of the solution demand prohibitively small time steps for explicit solvers. This work addresses the stiffness issue when employing neural ODEs for model order reduction by introducing a suitable reparametrization in time. The considered map is data-driven and it is induced by the adaptive time-stepping of an implicit solver on a reference solution. We show that the map produces a nonstiff system that can be cheaply solved with an explicit time integration scheme. The original, stiff, time dynamic is recovered by means of a map learnt by a neural network that connects the state space to the time reparametrization. We validate our method through extensive experiments, demonstrating improvements in efficiency for the neural ODE inference while maintaining robustness and accuracy when compared to an implicit solver applied to the stiff system with the original right-hand side.
\end{changemargin}
}
\vspace{0.2cm}
\noindent\textbf{Key words:} Ordinary differential equations, Stiff equation, Reduced order model, Neural networks, Runge-Kutta. 

\noindent\textbf{AMS subject classification:} 65L99, 68T07.


\section{Introduction}
In computational science and engineering, the complexity and scale of numerical simulations have grown exponentially, driven by advances in hardware and algorithms. High-fidelity models {-- also called Full Order Models (FOMs) --} can now capture intricate details of physical phenomena across diverse domains such as fluid dynamics, structural analysis, climate modeling, and materials science. However, the increasing complexity of these models often comes at a significant computational cost, rendering them impractical for many applications that require repeated model evaluations over a large number of parameter values. Indeed, the computational cost is driven up by the curse of dimensionality \cite{benner2015survey}.

Reduced Order Models (ROMs) have emerged as a powerful solution to this challenge, offering a balance between accuracy and computational efficiency. By compressing the essential features of {FOMs} into a more manageable form, usually called latent space, ROMs enable the efficient simulation and analysis of complex systems, especially those requiring real-time analysis \cite{peirlinck2021precision}, inverse problems \cite{cui2015data, frangos2010surrogate, maday2015parameterized, lahivaara2019estimation}, uncertainty quantification \cite{bui2008parametric, sudret2000stochastic, chen2017reduced}, shape optimization \cite{manzoni2012shape} or optimal control \cite{negri2013reduced, ravindran2000reduced, bergmann20087813}.

The field of ROMs has seen substantial growth, providing practitioners with a wide range of techniques to choose from. Each method comes with its own set of advantages and drawbacks, and the optimal choice of a ROM often depends on the specific requirements of the application. {Indeed, the term ROMs refers to multiple approaches, which are rather different from one another.}

One of the most well-established paradigms for model order reduction is proper orthogonal decomposition (POD) \cite{hesthaven2016certified, benner2017model}, that is, a linear projection technique that provides a high compression rate and, for instance, achieves a particularly high level of accuracy in the case of diffusion processes.

If the governing equations of the high-fidelity model are explicitly employed in the ROM, they are called intrusive \cite{prud2002reliable, hesthaven2022reduced, benner2005dimension, antoulas2005approximation, bui2008model}. Usually, a greedy algorithm or POD is employed to project the {FOM} onto a lower-dimensional space. However, for nonlinear models, additional techniques like the (discrete) empirical interpolation method are required to handle the complexities, which can lead to an unfavorable trade-off between accuracy and computational efficiency \cite{barrault2004empirical, canuto2009posteriori, chaturantabut2010nonlinear}. Moreover, in \cite{lee2020model} it has been shown that for transport-dominated problems, a large latent dimension must be used to achieve reasonable accuracy when using POD, thus limiting their practical use.

In recent years, data-driven methods based on machine learning and deep learning gained large traction due to advances in GPU-based hardware and developments in algorithms and software supporting the artificial intelligence ecosystem. This paradigm is based on the availability of a large quantity of high-fidelity data that can be used to train the model. This approach proves effective if the large amount of resources employed in the offline phase (generation of the solutions and training of the model) is repaid by a large enough number of (cheap) evaluations of the trained model (online phase). These models are called non-intrusive since only knowledge of the data is needed to build them \cite{hesthaven2018non, wang2019non}. {We also refer to these kinds of ROMs as surrogate models, given their ability to approximate the behavior of the original FOM from the data.} The field is ever-growing, with a large number of techniques being developed and published. To the best of our knowledge, the most relevant examples of data-driven approaches are the following.
Closure models are a hybrid between POD and data-driven approaches where a non-linear term (a neural network) is added to the linear ROM (usually POD-Galerkin) to account for the effect of the unresolved dynamics \cite{wang2020recurrent}.
In \cite{oommen2022learning,vlachas2022multiscale,brunton2020machine,maulik2021reduced,fresca2021comprehensive}, densely connected or convolutional autoencoders have been used to extract the full latent variables. Other techniques have been proposed to enhance the accuracy. For instance, graph neural networks (GNNs) tackle the problem of complicated geometries discretized with unstructured meshes, offering a more natural representation of variables using a graph \cite{pichi2024graph, pegolotti2024learning}. Fourier neural operators are a very successful approach that instead exploits the frequency space to enrich the layers of the model \cite{li2020fourier, li2023fourier, wen2022u}. 
Depending on the considered model, time evolution may be treated just as any parameter of the system or explicitly learned, examples are DeepONets \cite{lu2021learning, oommen2022learning}, dynamic mode decomposition \cite{schmid2010dynamic, brunton2022data}, sparse identification of reduced latent dynamics (SINDy) \cite{champion2019data, brunton2016discovering}, RNNs \cite{liu2022hierarchical}, LSTMs \cite{vlachas2022multiscale, conti2023multi} or Gaussian processes \cite{guo201975}.

Another time evolution approach of particular interest is based on neural ODEs \cite{chen2018neural, regazzoni2019machine, Regazzoni2024, portwood_turbulence_2019}. The core concept of this method is straightforward: it involves learning the right-hand side of a dynamical system with a dense neural network, which can be viewed as the continuous limit of a residual neural network. This approach offers numerous advantages, including a more natural representation of continuous-time data, and enhanced parameter and memory efficiency.
The primary challenge with neural ODEs is that they are in all effects ODEs that must be solved. This often leads to issues with stiffness, complicating their solution and potentially affecting their stability and accuracy. Stiffness is the condition where there are rapid changes in some components of the solution while other components change slowly. This phenomenon creates numerical challenges when solving the equations, as explicit numerical methods might become inefficient or unstable \cite{wanner1996solving}. Thus, the necessity of employing implicit solvers significantly drives up the cost.

In the present work, we tackle the problem of building a ROM of a stiff system. In this case, the main bottleneck is not represented by the large number of variables of the system but by the stiffness of the system.  Indeed, {there are multiple types of} stiffness, which are usually connected to features of the system, such as the presence of multiple time scales, the spectrum of the Jacobian of the system, or numerical instability, that make building a ROM a challenging task. 
In this case, employing fixed-time timesteppers is unfeasible since the rapidly varying scales require unreasonably small steps. On the other hand, treating time like an input of a dense neural network also proves to be problematic since neural networks are known to perform badly on high frequency \cite{WANG2022110768}. Still, in literature, there can be found several examples of stiff problems solved with machine learning. A possible solution was first proposed in \cite{ji2021stiff}, where the authors employ the quasi-steady-state assumption (QSSA) to reduce the stiffness of the ODE systems and show that a physics-informed neural network (PINN) then can be successfully applied to the converted non-/mild-stiff systems. This is particularly interesting since the system can be solved with an explicit method. In \cite{kim2021stiff}, the authors show that computing a stabilized gradient and suitable scaling of the network outputs enables learning stiff neural ODEs. In \cite{deflorio2022physics}, PINNs are used with extreme learning machine to solve stiff problems. In \cite{ANANTHARAMAN20221}, stiff quantitative systems pharmacology models are accelerated with echo state networks. In \cite{Zhang_2022}, an autoencoder is used to produce a reduced order model of a chemical kinetic model for the simulation of a combustion system. Recently, in \cite{GOSWAMI2024116674} DeepONets are used to learn the discretized solution of challenging stiff chemical kinetics. 

{Our strategy is to employ a neural ODE-based non-stiff surrogate for stiff dynamics. Namely, we} aim to employ a neural ODE to learn a time-reparametrized stiff system, where the time map is suitably built to reduce the stiffness of the problem. The insight is to exploit the advantages of the neural ODEs we previously described while keeping the (online) cost of solving the ODE under control by using an explicit solver. Indeed, our main contribution is to employ the time reparametrization induced by the time-stepping of an implicit solver as a way to reduce the stiffness of the system. An advantage of this procedure is that it is purely data-driven, meaning that it is possible to avoid manually deriving explicit algebraic expressions for the QSS species as it was done in \cite{ji2021stiff}. On the other hand, a drawback of our approach is that it is only worthwhile for ROMs: the expensive phase of gathering data for building the time mapping must be offset by a cheap online phase. Our technique can also be interpreted as a way to embed adaptive (implicit) time-stepping directly into the neural ODE. Thus, we can avoid backpropagation through the implicit solver, which incurs a cubic cost with respect to the number of neural network parameters \cite{kim2021stiff}, thereby greatly reducing the computational cost of training the neural ODE. Our reparametrization strategy is also connected with \cite{van2023adaptive,barral2024registration,khamlich2024optimal}, where the authors leverage spatial mappings for adaptive mesh refining to improve computational efficiency and to enhance the representation of nonlinear dynamics in complex systems.
{Moreover, recent work \cite{loya2025structure} has demonstrated that a structure-preserving NODE with linear/nonlinear splitting can achieve long-term stability on stiff problems comparable to implicit methods. 
}
Once the nonstiff system is solved, the solution is mapped back to the stiff dynamics using a map learnt with a neural network. The neural network is trained once with the neural ODE and shows good generalization properties to parameters outside the training dataset. Particular care is taken in the definition of the map so that it has no dependence on time. Hence, the model maintains good generalization in time for periodic systems since there is only explicit dependence on the state space.

{Our technique is designed especially for problems where stiffness can be traced back to the case of large eigenvalues with a negative real part. Indeed, another largely studied type of stiffness is high-frequency stiffness, which corresponds to the case of reasonably large eigenvalues but very high frequency \cite{petzold1997numerical}. In these cases, the system features a fast and a slow solution, where the fast solution must be solved accurately to track the slow one. While reparametrization approaches have been proposed \cite{huang1997adaptive}, they have significant limitations since they need to resolve the fast frequency \cite{ascher1999some, richardson2011symplectic}. Thus, this kind of stiffness remains a challenge for our approach and an area of ongoing research. Indeed, other authors have explored structure-preserving machine learning techniques, such as fast-slow neural networks, to bypass the direct resolution of fast transients in singularly perturbed systems \cite{serino2024intelligent}.
}

We demonstrate the accuracy and speed of the proposed ROM on a set of five test problems, where explicit solvers usually fail, widely used in literature to benchmark stiff ODE solvers \cite{mazzia2003test}. The results are compared with a state-of-the-art Runge-Kutta implicit solver of the Radau II A kind of the fifth order \cite{wanner1996solving} applied to the stiff system with the original right-hand side. 

\vspace{3mm}
The remainder of the paper is organized as follows. First, in Section~\ref{sec:problem formulation} we introduce the mathematical settings and notation for the problems of interest. In particular, in Section~\ref{sec:neural-odes} we review the neural ODEs approach, which {is} the foundation of this work. In Section~\ref{sec:stiffness} we introduce our methodology, discussing also the implementation details. In Section~\ref{sec:numerical-results} we present a comprehensive set of numerical results. Finally, in Section~\ref{sec:conclusions} we draw some conclusions, discussing the strengths and limitations of the proposed approach.


\section{Problem formulation}\label{sec:problem formulation}
In this section, we first present the construction of ROMs in algebraic terms for approximating a parametric dynamical system. Then, we provide a concise overview of fundamental deep learning concepts and formulate the data-driven model order reduction problem within the framework of neural ODEs. 

\subsection{Parametric systems of ordinary differential equations}
We consider the following parameterized, finite-dimensional dynamical system, described by a set of first-order ODE:

\begin{equation}
    \begin{cases}
        \dot{\mbf u}(t; \mbf \mu) = \mbf f (t, \mbf u; \mbf \mu), \qquad t \in (0, T]\\
        \mbf u(0; \mbf \mu) = \mbf u_0(\mbf \mu).
    \end{cases}
    \label{eq:dynamic-parametric-system}
\end{equation}
where $\mbf \mu$ is a vector containing all the parameters of the system, which belongs to a compact set $\Gamma \subset \mathbb R^{N_\mu}$, $\mbf u: [0, T] \times \Gamma \rightarrow \mathbb R^{N_u}$ is the parameterized solution, $\mbf f: \mathbb [0, T] \times \mathbb R^{N_u} \times \Gamma \to \mathbb R^{N_u}$ is the right-hand side function (encoding the dynamical system), $t \in [0, T]$ is the time variable, $\mbf u_0: \Gamma \rightarrow \mathbb R^{N_u}$ is the initial condition, and $\dot{\mbf u}$ is the total derivative with respect to time $t$. Without loss of generality, we are considering that the initial time for the system is zero; indeed there always exists a translation in time that brings us to this case.

This formulation encompasses not only ODEs but also problems stemming from the semi-discrete formulation of systems of partial differential equations (PDEs) through suitable methods (e.g.: finite elements \cite{zienkiewicz2005finite}, spectral elements \cite{canuto2007spectral}, discontinuous Galerkin \cite{cockburn2012discontinuous, hesthaven2007nodal}, etc.). In these cases, one may observe an increase in the computational cost of the simulation since the dimension of the state space $N_u$ is usually large. Indeed, these methods usually rely on a fine mesh to discretize the spatial derivatives. 

The meaning of the parameter $\mbf \mu \in \Gamma$ is broad, indeed it may represent physical properties (e.g.: material properties), initial conditions, or geometrical properties (e.g.: the shape of the domain). In this work, we focus on cases where $\mbf \mu$ models physical properties.

\subsection{Solving ordinary differential equations}\label{sec:runge-kutta}
There is a wide literature of methods to solve ODEs \cite{wanner1996solving}. In this work, we concentrate on the well-known Runge-Kutta methods. Let $r \in \mathbb{N}$ s.t.\ $r \geq 1$, let $b_i, a_{ij}, (i, j = 1, ..., r)$ be real numbers and let $c_i$ be defined as
\begin{equation*}
    c_i = \sum_{j=1}^{i-1} a_{ij}, \quad \forall \, i = 1, ..., r.
\end{equation*}
Then, fixed a timestep $\Delta t > 0$ and setting $t^{(n)} = n \Delta t$, we call the $r$-stage Runge-Kutta method the following discretization of Eq.~\bref{eq:dynamic-parametric-system}
\begin{equation}
    \begin{split}
            & \mbf k_i = \mbf f\left(t^{(n)} + c_i \Delta t, \mbf u^{(n)} + \Delta t \sum_{j=1}^r a_{ij} \mbf k_j\right), \quad i = 1, ..., r\\
            & \mbf u^{(n+1)} = \mbf u^{(n)} + \Delta t \sum_{i=1}^r b_i \mbf k_i,
    \end{split}
    \label{eq:rk}
\end{equation}
$\forall \, n \in \mathbb{N}$ s.t.\ $t^{(n)} \leq T$. Here, $\mbf u^{(n+1)} \approx \mbf u(t^{(n)} + \Delta t)$. If $a_{ij} = 0$ for $i \leq j$, each stage of Eq.~\bref{eq:rk} depends only on the previous ones, thus the method is explicit. If $a_{ij} = 0$ for $i < j$ and at least one $a_{ii} = 0$, the method is called diagonally implicit. In all the other cases, it is called implicit.

Of particular interest are embedded Runge-Kutta methods, which automatically adjust the step size $\Delta t$ to achieve a prescribed tolerance of the local error. The embedded method exploits a second set of scalars $b_i^*, i = 1, ..., r$ to compute another approximation of the solution of lower order $\mbf u^*$. Namely, we want to satisfy the component-wise inequality
\begin{equation*}
    | (\mbf u^{(n+1)})_i - (\mbf u^{*(n+1)})_i| \leq \delta_i, \qquad \delta_i = \texttt{atol} + \texttt{rtol}\max(| (\mbf u^{(n+1)})_i|,| (\mbf u^{(n)})_i|), \qquad i = 1, ..., N_u
\end{equation*}
where \texttt{atol} and \texttt{rtol} are user-defined (positive) absolute and relative tolerances, respectively. Then, the error is estimated as 
\begin{equation*}
    \varepsilon = \sqrt{\frac{1}{N_u} \sum_{i=1}^{N_u}\left( \frac{ (\mbf u^{(n+1)})_i -  (\mbf u^{*(n+1)})_i}{\delta_i} \right)^2}
\end{equation*}
and it is compared to one to find the next step size, namely
\begin{equation*}
    \Delta t^{(n+1)} = \Delta t^{(n)} \min(10, \max(0.1, \varepsilon^{-1 / (q + 1)})),
\end{equation*}
where $q$ is usually chosen to be the order of the method, and 10 and 0.1 are the maximum and minimum change factors, respectively. The value of these factors can be tuned depending on the problem. If the normalized error is larger than one, the timestep is rejected and is computed again with the updated $\Delta t$.

Explicit Runge-Kutta methods are generally unsuitable for the solution of stiff equations because their region of absolute stability is small \cite{wanner1996solving}. Indeed, in the present work, we will employ three different Runge-Kutta methods depending on the application. Namely, we use an implicit Runge-Kutta of the Radau II A kind of the fifth order \cite{HAIRER199993} for solving the stiff problems, and we use a Runge-Kutta of the fourth order (in the adaptive and fixed step size variants) for nonstiff problems. We employ the following table for the Radau II A method
\begin{equation*}
    \arraycolsep=1.1pt\def\arraystretch{1.5}
    \begin{array}{c|ccc}
        {\frac {2}{5}}-{\frac {\sqrt {6}}{10}}&{\frac {11}{45}}-{\frac {7{\sqrt {6}}}{360}}&{\frac {37}{225}}-{\frac {169{\sqrt {6}}}{1800}}&-{\frac {2}{225}}+{\frac {\sqrt {6}}{75}}\\
        {\frac {2}{5}}+{\frac {\sqrt {6}}{10}}&{\frac {37}{225}}+{\frac {169{\sqrt {6}}}{1800}}&{\frac {11}{45}}+{\frac {7{\sqrt {6}}}{360}}&-{\frac {2}{225}}-{\frac {\sqrt {6}}{75}}\\
        1&{\frac {4}{9}}-{\frac {\sqrt {6}}{36}}&{\frac {4}{9}}+{\frac {\sqrt {6}}{36}}&{\frac {1}{9}}\\
        \hline &{\frac {4}{9}}-{\frac {\sqrt {6}}{36}}&{\frac {4}{9}}+{\frac {\sqrt {6}}{36}}&{\frac {1}{9}}\\
    \end{array}  
\end{equation*}
Due to the super-convergence of the Radau II A method (classical order $p = 2r-1$) it is not possible to have an embedded method of order $p - 1$ without extra cost. By taking a linear combination of $\Delta t \mbf f(\mbf x_0, \mbf u_0)$ and the internal stage values $\mbf k_i$ it is however possible to get an approximation of order $r$, we refer the interested reader to \cite{HAIRER199993} for details. The following tableau is used for explicit Runge-Kutta of the fourth-order \cite{dormand1980family}, also called the Dormand-Prince method or RK45
\begin{equation*}
    \begin{array}{c|ccccccc}
        0&&&&&&&\\
        1/5&1/5&&&&&&\\
        3/10 &3/40 &9/40&&&&&\\
        4/5 &44/45 &-56/15 &32/9&&&&\\
        8/9 &19372/6561 &-25360/2187 &64448/6561 &-212/729&&&\\
        1 &9017/3168 &-355/33 &46732/5247 &49/176 &-5103/18656&&\\
        1 &35/384 &0 &500/1113 &125/192 &-2187/6784 &11/84&\\
        \hline
        &35/384&0&500/1113&125/192&-2187/6784&11/84&0\\
        &5179/57600 &0 &7571/16695 &393/640 &-92097/339200 &187/2100 &1/40 
    \end{array}  
\end{equation*}

\subsection{Data-driven model order reduction}
\begin{figure}[!t] 
    \centering
    \includegraphics[width=0.6\textwidth]{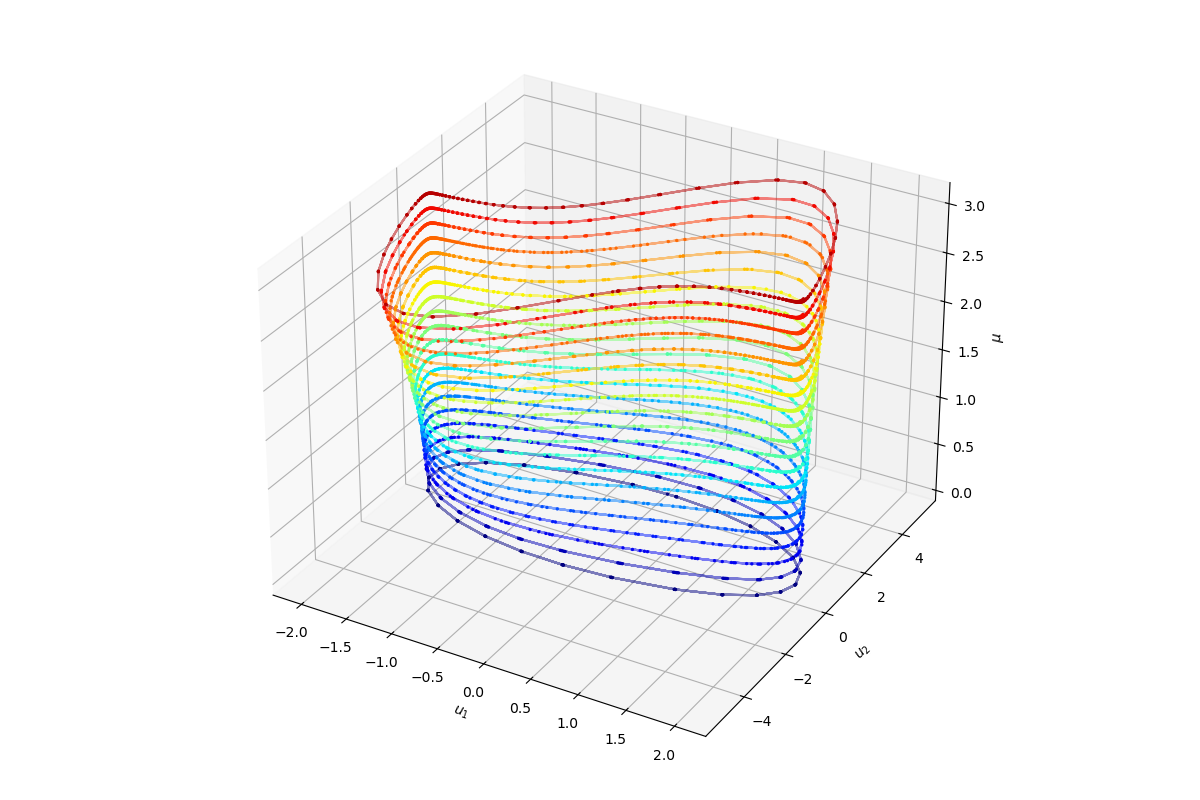}
    \caption{Example of discrete solution manifold $\widetilde{\mathcal{M}}$ for the Van der Pol oscillator.}
    \label{fig:maninfold-example}
\end{figure}

In this section, we discuss the concepts and notation that stand at the basis of data-driven model order reduction. We introduce the manifold of solutions of equation \bref{eq:dynamic-parametric-system}
\begin{equation}
    \mathcal{M} = \{\mbf u(t;\mbf \mu) \, | \, t \in [0, T], \mbf \mu \in \Gamma \} \subset \mathbb{R}^{N_u},
\end{equation}
which is a collection of the values of the solutions for each parameter $\mbf \mu \in \Gamma$ and time $t \in [0, T]$. We suppose that it is possible to find a discretized approximation $\mbf u(t^{(n)};\mbf \mu)$ for some times $t^{(n)} \in [0, T]$ up to arbitrary accuracy by integrating \bref{eq:dynamic-parametric-system} with a suitable method. By solving the problem for $N_\mu$ different choices of $\mbf \mu$ we create a discretized manifold

\begin{equation}
    \widetilde{\mathcal{M}} = \{
    \mbf u(t^{(n)}; \mbf \mu^{(j)}) \, | \, 
    \mbf \mu^{(j)} \in \Gamma, \, 1 \leq j \leq N_\mu, \,
    t^{(n)} \in [0, T], \, 1 \leq n \leq N(\mbf \mu^{(j)})
\},
\end{equation}
which is a finite set of vectors that approximate the solution at a certain discretized point in time $t^{(n)}$ and for a certain choice of parameter $\mbf \mu^{(j)}$. Let us stress that the time discretization $\{0 \leq t_1 < ... < t_{N(\mbf \mu^{(j)})} \leq T\}$ may change depending on the value of $\mbf \mu^{(j)}$, and need not be uniform. Indeed, it often happens that $\mbf \mu^{(j)}$ changes the dynamics of the system, and thus a different time discretization is generated by the adaptive procedure. For the sake of simplicity, we omit the explicit dependence on $\mbf \mu^{(j)}$ of the time discretization. 

In practice, one might obtain data also from knowing the exact analytical solution of the problem for some specific values of $\mbf \mu$ and $t$ or from real-world measurements.
Data-driven ROMs aim to leverage the information contained in the discretized solution manifold $\widetilde{\mathcal{M}}$ to obtain in a fast and reliable way the solution $\mbf u(t^{(n)}; \mbf \mu^{(j)})$ for $\mbf \mu^{(j)}$ or $t^{(n)}$ not in the discrete manifold used for training the model. Figure~\ref{fig:maninfold-example} shows a graphical representation of the solution manifold for the Van der Pol oscillator (cf. Section~\ref{sec:van-der-pol}).
The motivation behind ROMs is that classical time integration methods may become prohibitively expensive either due to the size of $N_u$ or due to stiffness (cf. Section~\ref{sec:stiffness}).

\subsection{Neural ODEs}\label{sec:neural-odes}
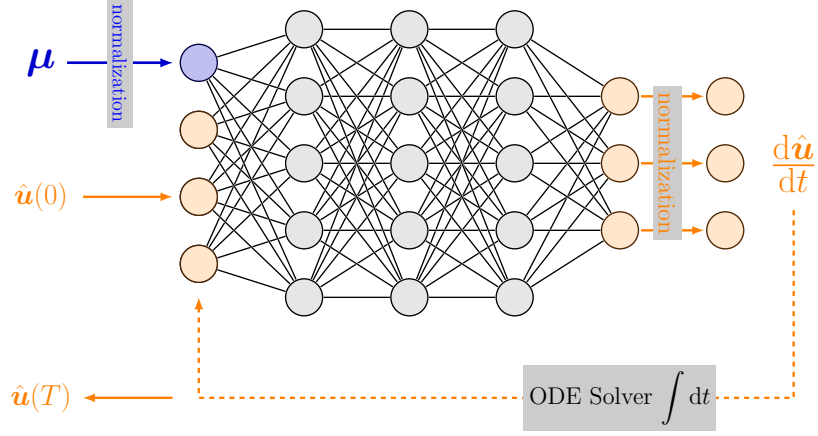
\begin{figure}[!t] 
    \centering
    \resizebox{0.7\linewidth}{!}{\input{tikz/fnn.tex}}
    \caption{Architecture of the Neural ODE for model order reduction. }
    \label{fig:neural-ode}
\end{figure}

The model that we choose to use to build the ROM is neural ODEs \cite{chen2018neural}. Neural ODEs have emerged as powerful tools to model continuous-time dynamics using a differential equation. This approach provides several key benefits. Firstly, neural ODEs offer superior flexibility in modeling time-series data, as they can naturally handle irregularly sampled data and varying time steps. Secondly, they ensure smooth trajectories and better integration with physical laws and other domains where the underlying processes are continuous in nature. Furthermore, neural ODEs can lead to more efficient training and inference, as they require fewer parameters and can leverage advanced numerical solvers for differential equations. Neural ODEs are also more memory efficient compared to traditional recurrent neural networks, especially for long time-series data, because they avoid the need to store intermediate states of the network during backpropagation.

More precisely, a neural ODE is {a dense fully connected feedforward (FNN) neural network $\mathcal{NN}$ \cite{goodfellow2016deeplearning,zhang2023dive}}, with parameters $\mbf \theta$, that evolves the dynamics as 
\begin{equation}
    \dot{\mbf u}(t; \mbf \mu) = \mathcal{NN}(t, \mbf u(t); \mbf \mu),
\end{equation}
that is, the FNN is employed to approximate the right-hand side $\mbf{f}$ of the dynamical system \bref{eq:dynamic-parametric-system}, thus it takes as inputs the time $t$, the current state $\mbf u$, and the parameters $\mbf \mu$. In other words, this technique aims to predict $\mbf u$ by using $\mathcal{NN}$ to learn the $\mbf{f}$ of Eq.~\bref{eq:dynamic-parametric-system}.
One of the main drawbacks of this approach is that it is necessary to numerically integrate the system in the online phase. Given the recurrent nature of the architecture, adjoint solvers have been developed to enhance the memory efficiency during the training phase \cite{chen2018neural}; their cost however remains much larger when compared to FNN.

When used as an architecture for ROMs, the neural ODE is often paired with another neural network to reduce the dimension of the system from $N_u$ to a small latent state space. This is usually achieved by means of autoencoders \cite{oommen2022learning,vlachas2022multiscale,brunton2020machine,maulik2021reduced,fresca2021comprehensive} or other reduction techniques \cite{Regazzoni2024}. However, if $N_u$ is small the advantages of using a small latent space representation may not be counterbalanced by the increased cost of evaluating the neural network. In this work, we propose a solution to this latter problem by introducing an appropriate change of variables that significantly reduces the stiffness of the problem. Indeed, we will show that our method makes it possible to use explicit solvers to integrate neural ODEs representing stiff systems.

We report in Figure~\ref{fig:neural-ode} a representation of the neural ODE. The hat symbol $\hat{\cdot}$ represents a suitable normalization of the data, see the next section for details. For the sake of simplicity, we will indicate in the captions the reference solution (the one that is comprised in the dataset) with the subscript ``ref'', and the components of $\mbf u$ as predicted by the neural ODE with $u_i$.

\subsubsection{Training a Neural ODE}

\begin{figure}[t]
    \centering
    \begin{minipage}{0.5\textwidth}
        \centering
        \hspace{0.7cm} Supervised training
        \includegraphics[trim={0 0.9cm 0 0.8cm},clip,width=0.7\linewidth]{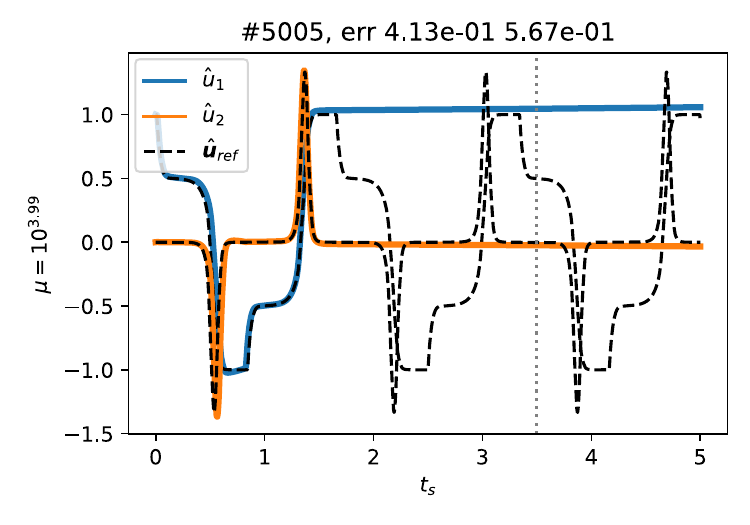}
    \end{minipage}%
    \begin{minipage}{0.5\textwidth}
        \centering
        \hspace{0.7cm} ODE solver training
        \includegraphics[trim={0 0.9cm 0 0.8cm},clip,width=0.7\linewidth]{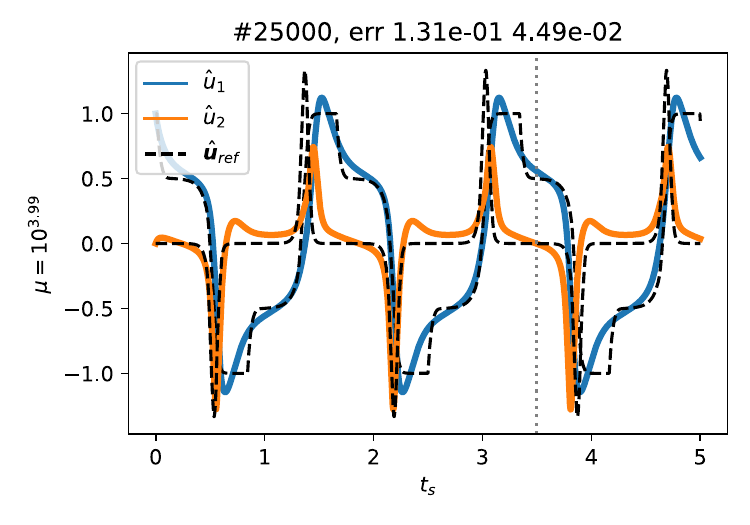}
    \end{minipage}
    \caption{Comparison of the supervised and ODE solver approaches for training the neural ODE from the same randomly initialized model for the Van der Pol oscillator (cf.\ Section~\ref{sec:van-der-pol}). In blue and orange the predictions of the model, and in black the reference solution. The supervised training is globally less accurate and stable, but better captures sharp gradients.}
    \label{fig:supervised-vs-node}
\end{figure}

\begin{figure}[t]
    \centering
    \begin{minipage}{0.5\textwidth}
        \centering
        Supervised training
        \includegraphics[trim={0 0 0 0},clip,width=0.9\linewidth]{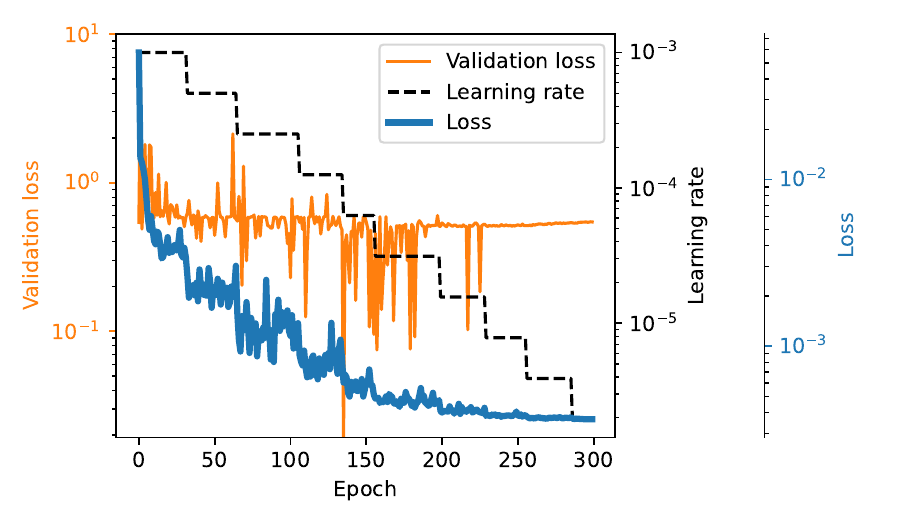}
    \end{minipage}%
    \begin{minipage}{0.5\textwidth}
        \centering
        ODE solver training
        \includegraphics[trim={0 0 0 0},clip,width=0.9\linewidth]{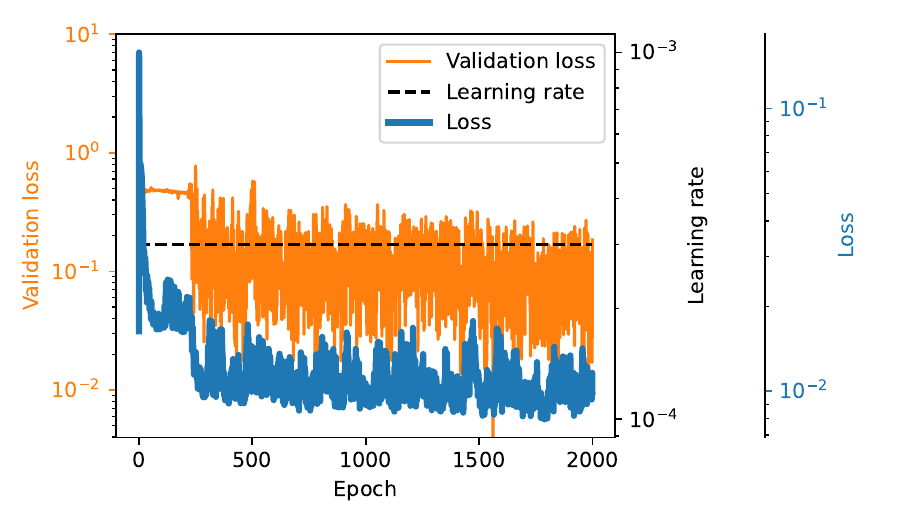}
    \end{minipage}
    \caption{Example of the behavior of the loss, the validation loss, and the learning rate during the training of the ROM with the two-step approach: first the supervised approach and then the ODE solver approach.}
    \label{fig:loss}
\end{figure}

There are two main approaches to train a neural ODE. The first one is the classic supervised approach. Assuming that the value of $\mbf f(t^{(n)}, \mbf u(t^{(n)}; \mbf \mu^{(j)}); \mbf \mu^{(j)})$ is known for each sample in the discrete manifold $\widetilde{\mathcal{M}}$, then it is possible to compute the MSE loss and minimize it through a gradient-based optimizer. However, this kind of data is not always available. Instead, the second approach -- the one originally proposed for the training of neural ODEs -- minimizes the error with respect to the timeseries data $\mbf u(t^{(n)}; \mbf \mu^{(j)})$ for a fixed $j$. Namely, by using a differentiable ODE solver, we minimize the loss

\begin{equation}
    \mathcal{L}(\mbf{\theta}) = \frac{1}{N_{\textnormal{train}}} \sum_{j = 1}^{N_{\textnormal{train}}} \norm{{\mbf{u}(t^{(N)}; \mbf \mu^{(j)})} - \left(\mbf{u}(t^{(0)}; \mbf \mu^{(j)}) + \int_{t^{(0)}}^{t^{(N)}} \mathcal{NN}(\tau, \mbf u(\tau; \mbf \mu^{(j)}), \mbf \mu^{(j)}; \mbf \theta) \mathrm{d}\tau \right)}_p,
    \label{eq:node-loss}
\end{equation}
where $\norm{\cdot}_p$ indicates the discrete $p$-norm, the elements $\mbf{u}(t^{(0)}; \mbf \mu^{(j)}) $ and $ {\mbf{u}(t^{(N)}; \mbf \mu^{(j)})}$ belong to the solution manifold $\widetilde{\mathcal{M}}$, $\mathcal{NN}$ is the approximation of $\mbf{f}$ made by a suitably trained neural network and the integral is usually computed with an explicit integration method like forward Euler or an explicit Runge-Kutta method. 

The advantage of the supervised approach is that it is less computationally demanding since it does not need to compute the gradient through the ODE solver using automatic differentiation. On the other hand, in the supervised framework, the loss is not always indicative of the accuracy of the predictions on the timeseries data $\mbf u(t^{(n)}; \mbf \mu), n=0, ..., N$. For this reason, we always employ as validation loss Eq.~\bref{eq:node-loss}. 

Implementation details for Eq.~\bref{eq:node-loss} are important when it is employed to train the neural network. We resample the discrete timeseries to match the output of the time integrator, namely a Runge-Kutta method of the fourth order with a fixed timestep. To minimize the interpolation errors, we use a reference solution with a very fine time discretization. A uniform resampling allows for an efficient use of the partial results of the ODE solver for all the times smaller than the final one. Moreover, the initial and final integration times are not $(t^{(0)}, t^{(N)})$ but are batched by uniform extraction between $[0, T]$ so that each time integration in the dataset is done with the same number of timesteps: this procedure enables vectorization of the computations, greatly reducing the computational cost. This optimization was originally introduced in \texttt{torchdiffeq} \cite{chen2018neural}. The number of timesteps used is called unroll length and is a hyperparameter that plays a key role in the optimization procedure. Finally, as commonly done in supervised learning, the computations are split into mini-batches to increase stochasticity.
If not otherwise stated, when using Eq.~\bref{eq:node-loss} for validation or testing we employ as a time integrator an embedded Runge-Kutta method of order four with adaptive time-stepping. However, computations are significantly cheaper since we can turn off automatic differentiation. 

In our experiments, we employ a two-phase training. First, a standard supervised approach is used to ``initialize'' the neural network since it is cheap and enables a fast reduction of the loss. {The derivative data is estimated with a filtering approach from the trajectory data (cf.\ Section~\ref{sec:numerical-results}).} Then, we switch to the differentiable ODE solver approach to fine-tune the ROM. Our experiments show that in this case, it is beneficial to use an increasingly larger unroll length to ensure the robustness and accuracy of the model. {The two-phase training is not strictly needed, but we have found it significantly reduces the time and effort required for hyperparameter tuning compared to using only the trajectory data.} Further details will be given in Section~\ref{sec:numerical-results} for each test case. Employing a two-phase training approach not only reduces the computational cost but also increases the quality of the final model. Indeed, from our experiments, we have observed that with the supervised learning approach the neural network is able to learn sharper gradients that may lead to unstable states. On the other hand, learning through the ODE solver produces stable results that might follow less precisely the reference solution. Figure~\ref{fig:supervised-vs-node} shows an example for two (underfitting) models. In Figure~\ref{fig:loss} we show the training and validation loss for a model trained in sequence with the supervised and ODE solver approaches.

Another important point of discussion is the explicit dependence on time $t$ of the $\mathcal{NN}$. In practice, we do not use $t$ as input of the neural network $\mathcal{NN}$. However, this does not necessarily limit us to autonomous systems of ODEs. Indeed, it is possible to use an input layer with memory, like an LSTM, to add temporal dependence. This approach is usually preferred to employing $t$ as input of the neural network since it improves generalization properties in time.

To accelerate the training {it is key to apply} suitable normalizations to the input $\mbf \mu$, the output $\mbf f$, and the timeseries data $\mbf u$. The exact expressions for the normalizations depend on the test case considered and are detailed in Section~\ref{sec:numerical-results}. {In some special cases, e.g., if $\mbf \mu$ is high-dimensional, further treatment of $\mbf \mu$ may be needed to accelerate the training.} The insight behind the choice of the normalization is that we aim to have data that has a magnitude close to one, does not have outliers, and has either a uniform distribution or a normal distribution with a standard deviation of magnitude close to one. In our experiments, strictly monotone normalizations achieved better results. We designed the normalizations by manual inspection of the histograms of the data.


\section{Addressing the stiffness issue} \label{sec:stiffness}
While many techniques have been developed to build reduced-order models of large systems stemming from semi-discretized PDEs, not many solutions are available for stiff ODEs. Indeed, it is a challenge to learn their complex dynamics, which often feature large gradients and widely different time scales. More rigorously, a commonly used formula to define the stiffness $S$ index of a problem is 

\begin{equation}
    S = \frac{\mathfrak{R}(\lambda_{\max{}})}{\mathfrak{R}(\lambda_{\min{}})}T,
\end{equation}
where $\lambda_i$ are the eigenvalues of the Jacobian of the system and $\mathfrak{R}$ indicates the real part.
However, it is known that this definition is not always helpful. For instance, it does not apply to the classic Robertson chemical kinetics problem since it has two zeros eigenvalues \cite{shampine2007stiff}. A more pragmatic definition is given by Harier: ``\textit{Stiff equations are problems for which explicit methods do not work}'' \cite{wanner1996solving}. That is, explicit numerical time integrators such as the embedded Runge-Kutta method of the fourth order require timesteps so small that implicit methods become cheaper. Stiffness is a particularly important issue when working with neural ODEs. Indeed, it is a well-known fact that neural ODE training may fail due to stiffness. For instance, \cite{kim2021stiff} shows issues of this kind and proposes techniques to improve the training. However, it is still necessary to use a costly implicit solver in the online phase, since the learnt dynamic is stiff.

Hence, to use the neural ODE as a ROM, we aim to reduce the stiffness of the system, allowing for the use of (cheap) explicit methods. To this end, we propose to change the dynamics of the problem by considering the change of variable in time induced by the adaptive time-stepping of the implicit solver used to solve the 
{FOM. We remark that this change of variable is particularly effective in the case of stiffness characterized by a large $S$. However, in cases where stiffness is characterized by the presence of an oscillating fast solution, reparametrization techniques show several limitations \cite{ascher1999some, richardson2011symplectic}. Our approach is not an exception, and this is an area for further improvement.}

We suppose the existence of a change of variable $t_s = t_s(t)$ that makes the system less stiff. Namely, instead of learning Eq.~\bref{eq:dynamic-parametric-system}, we learn the following system 
\begin{equation}
 \begin{cases}
 \dot{\mbf u_s}(t_s; \mbf \mu) = \mbf f_s (t_s, \mbf u_s; \mbf \mu), \qquad t_s \in (0, 1]\\
 \mbf u_s(0; \mbf \mu) = \mbf u_{s,0}(\mbf \mu),
 \end{cases}
    \label{eq:reparam-ode}
\end{equation}
where, $\mbf u_s(t_s; \mbf \mu) = \mbf u(t(t_s); \mbf \mu)$ and thus $\mbf f_s$ is
\begin{equation*}
 \mbf f_s(t_s, \mbf u_s; \mbf \mu) = \frac{\text{d} \mbf u_s(t_s; \mbf \mu)}{\text{d} t_s} = \frac{\text{d} \mbf u(t(t_s); \mbf \mu)}{\text{d} t} \frac{\text{d} t}{\text{d} t_s} = \mbf f(t(t_s; \mbf \mu), \mbf u; \mbf \mu) \dot{t}(t_s; \mbf \mu).
\end{equation*}
We aim to build the time mapping $t_s = t_s(t)$ in a data-driven manner, in particular we suppose that
\begin{equation*}
 t_s(t^{(n)}) = \frac{n}{N},
\end{equation*}
that is, the $\Delta t$ chosen at the $n$-th step by the implicit time stepper is proportional to the derivative of the function that we employ as a change of variable in time. We remark that $t_s$ is strictly increasing and thus invertible, with inverse $t = t(t_s)$. We guess this change of variables induces less stiff dynamics since the adaptive time-stepping spreads the sharp peaks and quickly adapts to the large timescales. The insight behind this approach is that using an implicit solver to train the neural ODE on Eq.~\bref{eq:dynamic-parametric-system} has a cubic cost \cite{kim2021stiff}, instead, this change of variable allows the training via an explicit method, which is significantly cheaper. In other words, the adaptive timestepping is embedded into the neural ODE, avoiding the backpropagation step through the implicit solver. While the proposed approach is closely related to classical embedded methods, it reduces significantly the cost of training a neural ODE. 

While it could be possible to learn $t(t_s; \mbf \mu)$ and $\dot{t}(t_s; \mbf \mu)$ by simple regression methods (such as monotonicity-preserving interpolation), this approach has two significant drawbacks for our use case: (\textit{i}) it becomes unfeasible to use Eq.~\bref{eq:node-loss} to train the neural ODE; (\textit{ii}) these interpolation schemes do not generalize well for $t > T$ and quickly become expensive as the dimension of $\mbf \mu$ increases.

Thus, in this framework, we have two neural networks. The first, $\mathcal{NN}$ is used to approximate $\mbf f_s$, which only requires pointwise estimates of these values. Further details for the training are given in Section~\ref{sec:impl-details}. The second ingredient needed to accommodate for the time reparametrization is a way to
compute $t(t_s; \mbf \mu)$ in the online phase, to map the solution back to the original dynamics. Our experiments show that learning directly the map $t = t(t_s; \mbf \mu)$ is a hard task for a neural network, probably due to the presence of abrupt changes in the function. Instead, we found that learning $\dot{t}(\mbf u_{s}(t_s; \mbf \mu))$ is a stable approach that produces more accurate results. {We note that, provided $t = t(t_s; \mbf \mu)$ is continuously differentiable (or just absolutely continuous, if we only require the derivative almost everywhere), it can always be represented -- and hence learned via the neural ODE approximation -- as 
\begin{equation*}
    t(t_s; \mbf \mu) = \int_{t_s^{(0)}}^{t_s} \dot t(\tau; \mbf \mu) \text{d}\tau \approx \int_{t_s^{(0)}}^{t_s} \mathcal{NN}_t(\tau, \mbf u_s(\tau; \mbf \mu), \mbf \mu) \text{d}\tau,
\end{equation*}
by the fundamental theorem of calculus.
In this case, $\mbf u_s(\tau; \mbf \mu)$ is just an extra input used to accelerate the training. However, to improve the generalization properties and keep the system autonomous, we prefer to remove the explicit dependence on $t_s$. This is possible only if the trajectories of the ODE are an injective immersion in $\mathbb{R}^{N_u}$ or if the ``velocity'' $\dot{t}(t_s; \mbf \mu)$ is uniquely well defined at every point of the trajectory. Since it is non-trivial to prove that this second property holds, we have only empirically tested the hypothesis by analyzing the data. Indeed, in practice, this approach works even if the trajectories satisfy the previous properties almost everywhere. Hence, in all our test cases, }
the time mapping is computed with the following path integral
\begin{equation*}
    t(t_s; \mbf \mu) \approx \int_{\mbf u_s(t_s^{(0)})}^{\mbf u_s(t_s)} \mathcal{NN}_t(\mbf u_s; \mbf \mu) \text{d}\mbf u_s,
\end{equation*}
where the values of $\mbf u_s$ are the ones obtained by solving Eq.~\bref{eq:reparam-ode} by time integration of the neural ODE $\mathcal{NN}$.
Our experiments (see Section~\ref{sec:numerical-results}) demonstrate that the method generalizes well to parameters $\mbf \mu$ not present in the training set. The evaluation of the time dynamics is cheap since it can be done in one batch using vectorization. Indeed, in the first step, we compute $\mbf u_{s}(t_s; \mbf \mu)$, then, we carry out the cumulative integration of $\dot t(\mbf u_{s}(t_s; \mbf \mu))$ with Simpson's rule. In our experiments, this step requires only a few milliseconds and its cost is negligible with respect to the cost of the integration of the neural ODE. The main advantage of this approach is that it is cheap compared to autoencoders, especially for small problems since the only overhead is the computation of $t = t(t_s; \mbf \mu)$. However, it is also possible to combine the autoencoder with the described methodology: first, the autoencoder extracts a latent state representation, then, we could apply the time reparametrization in the latent space, reducing the stiffness of the latent dynamics. 

Unlike traditional non-constant time stepping, our data-driven time reparametrization fundamentally alters the underlying time dynamics by introducing a new, normalized time variable $t_s$ based on the adaptive time-stepping behavior of the original system. While the results are in practice similar, this technique enables to predict the adaptive timestepping for a new model parameter $\mbf \mu$ in a much cheaper way, without the necessity of evaluating the Jacobian or solving a linear system, like it is requested in an implicit ODE solver.

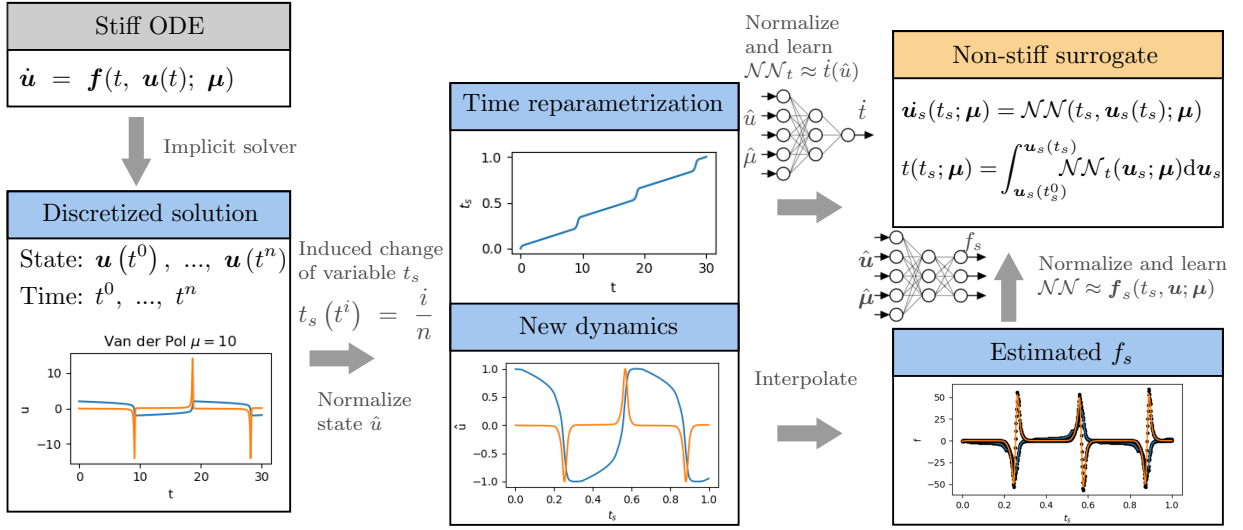
\begin{figure}[t] 
    \centering
    \resizebox{1.0\linewidth}{!}{\input{tikz/reparam_workflow.tex}}
    \caption{Time reparametrization workflow. Starting from the stiff ODE system (in gray), we compute the discrete solution and obtain the discrete manifold $\widetilde{\mathcal{M}}$ using an implicit method. We normalize the solution, compute the time reparametrization $t_s(t)$, and estimate $f_s$ (in blue). Finally, we are able to train a nonstiff model that can be solved using an explicit method (in orange).}
    \label{fig:workflow_dynamics_change}
\end{figure}

\subsection{Implementation details}\label{sec:impl-details}
Applying this methodology to the training with the differentiable ODE solver does not require any change. Indeed, the values of the timeseries $\mbf u_{s}(n / N; \mbf \mu) = \mbf u(t^{(n)}; \mbf \mu)$ are known, they only have a different time discretization, that we know pointwise. The main difference concerns the supervised learning approach. As explained before, this step allows us to ``initialize'' the neural network in a cheap way and requires us to estimate both $\dot{\mbf u}_{s}(t_s; \mbf \mu) = \mbf f_{s}(t_s, \mbf u_{s}; \mbf \mu)$ and $\dot{t}$, in order to train two neural network $\mathcal{NN}$ and $\mathcal{NN}_t$, respectively. The approximation of the derivative is done using a Savitzky-Golay filter with low order (two or three) and small window size ($<10$). Figure~\ref{fig:workflow_dynamics_change} shows a scheme of the workflow. To accelerate the training also the values of the time derivative are normalized, in particular, we have found a logarithmic normalization to work particularly well.

The optimizer we employ is AdamW \cite{loshchilov2017decoupled} with batch size 32 and a learning rate scheduler that geometrically reduces errors on plateaus of a factor of one-half with patience twenty epochs. The training lasts for a maximum of 400 epochs. The best model is chosen using early stopping based on the validation loss. The activation functions, the depth and width of the neural network, the initial learning rate, and the seed are hyperparameters that we optimize by means of random search on a problem-by-problem basis. In particular, the activation function is chosen among GELU \cite{hendrycks2016gaussian}, SiLU \cite{elfwing2018sigmoid}, Hard Swish \cite{howard2019searching}, Leaky ReLU and ReLU, the depth between three and ten, the width between five and 100 and the initial learning rate between $2\cdot10^{-3}$ and $10^{-4}$.

We employ a very similar approach to fine-tune the differentiable solver. The main difference is that the maximum number of epochs considered is 2000. The unroll length and the timestep size are tuned for each problem. As a rule of thumb we employ increasingly larger unroll lengths with smaller and smaller learning rates. Namely, we start with covering about 5\% of the time domain and it increases up to 50\%.

The neural network $\mathcal{NN}_t$ approximating $\dot t(\mbf u_{s}(t_s; \mbf \mu))$ is also trained in a two-step process. First, it is trained in a supervised manner. Similarly to what we do in the other supervised step, we employ as validation loss the accuracy of the predictions on the validation timeseries when employing the ODE solver to integrate both $\mathcal{NN}_t$, which approximates $\dot t(\mbf u_{s}(t_s; \mbf \mu))$, and $\mathcal{NN}$, which approximates $\mbf f_{s}(t_s, \mbf u_{s}; \mbf \mu)$. Then, it is fine-tuned using as a target the time $t$. This is achieved by automatic differentiation through Simpson's integration rule. Namely, the loss is
\begin{equation*}
    \mathcal{L}(\mbf{\theta}) = \frac{1}{N_{\textnormal{train}}} \sum_{j = 1}^{N_{\textnormal{train}}} \norm{t^{(N)} - \int_{\mbf u_{s}(t_s^{(0)})}^{\mbf u_{s}(t_s^{(N)})} \mathcal{NN}_t(\mbf u_{s}(\tau; \mbf \mu^{(j)}), \mbf \mu^{(j)}; \mbf \theta) \mathrm d \tau}_p,
\end{equation*}
where $t^{(0)}, t^{(N)}$ are sample times, and $\mbf{\mu}^{(j)}$ are sample parameters used to build the discrete manifold $\widetilde{\mathcal{M}}$ and the integral is computed cumulatively on the discretization induced by the explicit Runge-Kutta method applied to the trained $\mathcal{NN}$. The time integration interval $(t^{(0)}, t^{(N)})$ and the sequence of states $\mbf u_s(\tau; \mbf \mu^{(j)})$ needed for the integration can be picked either from the solution manifold or the values obtained by applying the trained $\mathcal{NN}$ to the training dataset. We tested $p=1,2,4$, the best results were achieved with $p=1$.

Finally, to achieve the best performances in the online phase (inference time), the model is optimized using the features of modern deep learning libraries, such as oneDNN graph fusion \cite{paszke2017automatic}.


\section{Numerical results}\label{sec:numerical-results}
Our methodology is tested in five numerical experiments featuring stiff problems widely used in the literature to benchmark implicit ODE integrators \cite{mazzia2003test,wanner1996solving}. We do not consider dynamical systems derived from semi-discrete PDEs since it would require the use of autoencoders to extract the latent dynamics. Indeed, we aim to assess the behavior of our method as a stand-alone improvement to neural ODEs. Introducing other variables in the test cases, such as the error of the autoencoders may be detrimental to this task. We are confident that given the wide range of test cases presented, our methodology could be extended to larger state spaces or semi-discrete PDEs without problems. Indeed, being able to outperform classical ODE solvers on small problems without the need to use a simplified latent space representation is a more difficult task.

Computations were carried out on a laptop with an AMD${}^\textnormal{\tiny\textregistered}$ Ryzen${}^\textnormal{\tiny\textregistered}$ 7 PRO 7840U. For the sake of reproducibility, we rely on the standard (CPU) implementations found in SciPy \cite{2020SciPy-NMeth} of the Runge-Kutta methods for solving the original stiff system with the original right-hand side. Indeed, the evaluation of the original right-hand side is much cheaper than the evaluation of a neural ODE that approximates it. All the test cases feature systems that cost orders of magnitude more to be solved with explicit solvers than with implicit ones. Moreover, to ensure a fair and rigorous comparison with the neural network, we executed both methods on a CPU. While leveraging a GPU would likely enable the neural network to significantly outperform the full-order method, our goal was to maintain a balanced evaluation framework. We remark that all the presented metrics, tables, and figures refer to the evaluation of $\mathcal{NN}$ and $\mathcal{NN}_t$ on a testing dataset, that is with values of the parameter $\mbf \mu$ that are unseen by the neural network.

Before proceeding with the presentation of the test cases, let us introduce some metrics of interest that we will use to assess the performance of the ROMs.
\begin{itemize}
    \item \textbf{time}: elapsed CPU time in seconds to solve the system up to the final time $T$.
    \item \textbf{\# fev}: number of evaluations of the right-hand side (either $\mbf f$ or $\mathcal{NN}$).
    \item \textbf{\# jev}: number of evaluations of the Jacobian of $\mbf f$ (used for the implicit solver).
    \item \textbf{\# lu}: number of LU decompositions (used for the implicit solver)
    \item \textbf{MSE${}_{t_s}$}: mean squared error of $\mathcal{NN}$ in the $t_s$ reparametrization.
    \item \textbf{MSE}: mean squared error in the $t$ parametrization.
    \item $L^2$: relative integral error computed with Simpson's rule.
    \item $d_{\text{peak}}$: average distance of the peaks (any sample whose two direct neighbors have a smaller amplitude) between the reference solution and the low fidelity solution (used for periodic solutions).
\end{itemize}
Finally, we indicate with \textbf{tol} the tolerance of the solver. If not otherwise specified, it is intended that we set both the absolute (\texttt{atol}) and relative (\texttt{rtol}) tolerances of the embedded solver to this value.

\subsection{Test case 1: Van der Pol oscillator}\label{sec:van-der-pol}

\begin{figure}[t]
    \centering
    \begin{minipage}{0.5\textwidth}
        \centering
        \includegraphics[trim={0 0.9cm 0 0.8cm},clip,width=0.7\linewidth]{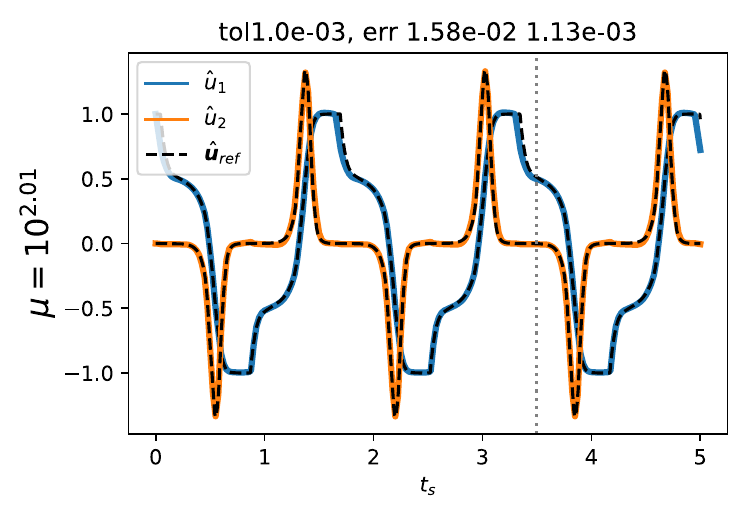}
    \end{minipage}%
    \begin{minipage}{0.5\textwidth}
        \centering
        \includegraphics[trim={0 0.9cm 0 0.8cm},clip,width=0.7\linewidth]{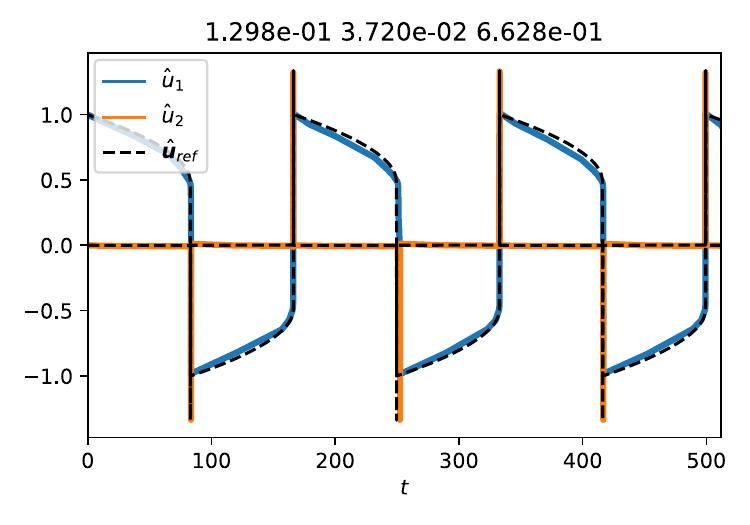}
    \end{minipage}
    \begin{minipage}{0.5\textwidth}
        \centering
        \includegraphics[trim={0 0.9cm 0 0.8cm},clip,width=0.7\linewidth]{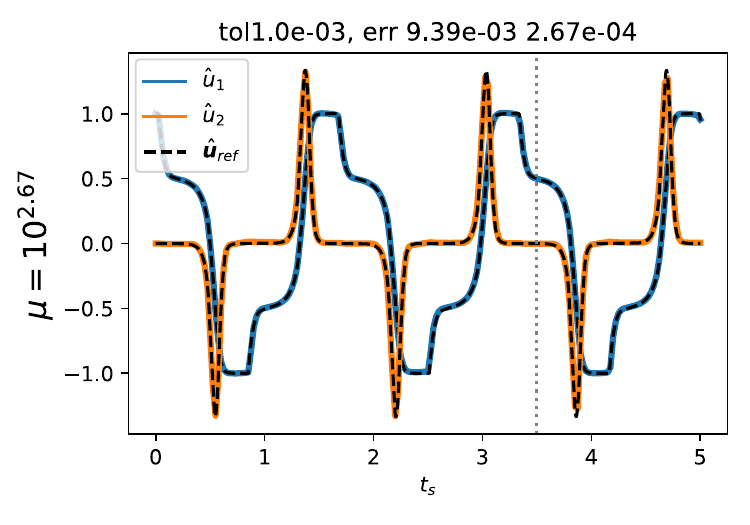}
    \end{minipage}%
    \begin{minipage}{0.5\textwidth}
        \centering
        \includegraphics[trim={0 0.9cm 0 0.8cm},clip,width=0.7\linewidth]{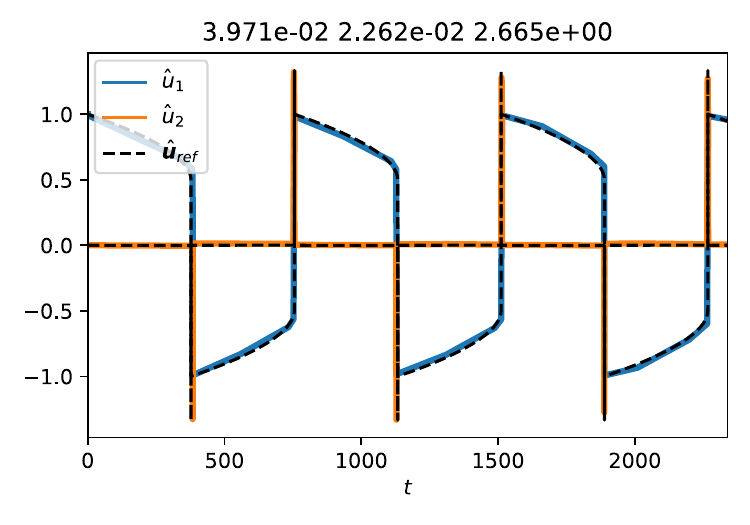}
    \end{minipage}
    \begin{minipage}{0.5\textwidth}
        \centering
        \includegraphics[trim={0 0.9cm 0 0.8cm},clip,width=0.7\linewidth]{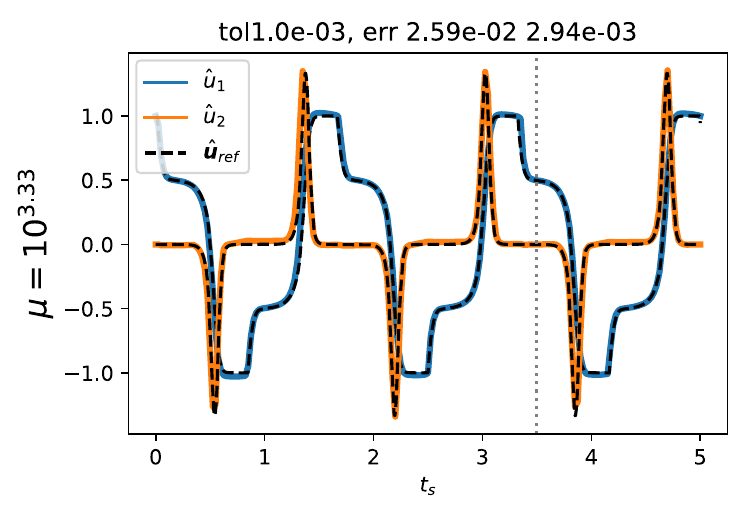}
    \end{minipage}%
    \begin{minipage}{0.5\textwidth}
        \centering
        \includegraphics[trim={0 0.9cm 0 0.8cm},clip,width=0.7\linewidth]{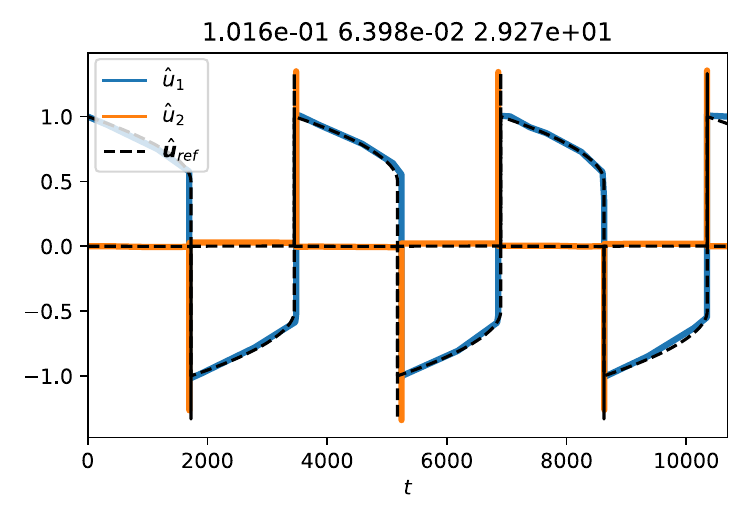}
    \end{minipage}
    \begin{minipage}{0.5\textwidth}
        \centering
        \includegraphics[trim={0 0 0 0.8cm},clip,width=0.7\linewidth]{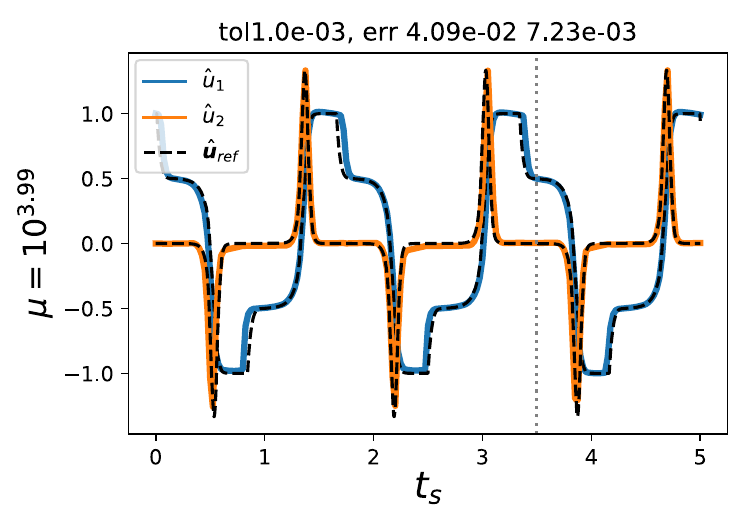}
    \end{minipage}%
    \begin{minipage}{0.5\textwidth}
        \centering
        \includegraphics[trim={0 0 0 0.8cm},clip,width=0.7\linewidth]{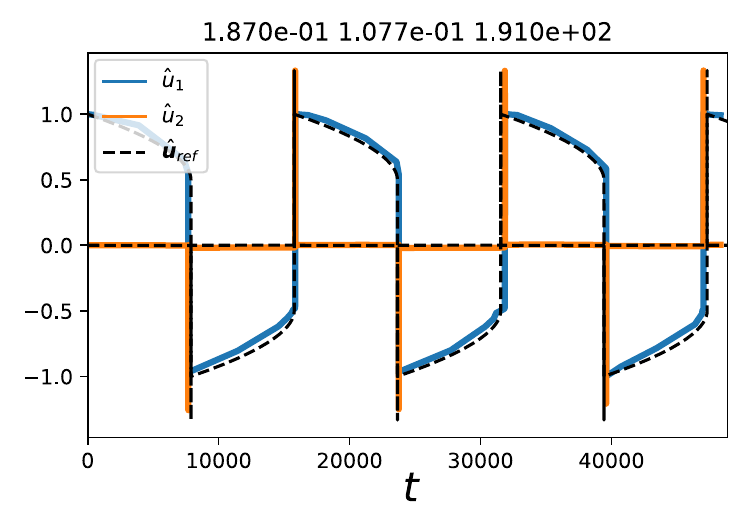}
    \end{minipage}

    \caption{Test case 1: Van der Pol oscillator. On the left, in colors, the neural ODE predictions in $t_s$. On the right, in colors, the prediction is mapped to the original time $t$. In black is the reference solution. The vertical line (\includegraphics[trim={0 0.5cm 0 0.5cm},width=0.5cm]{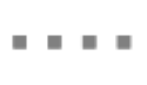}) represents the final time in the training dataset.}
    \label{fig:rom-vdp}
\end{figure}

The van der Pol oscillator is a second-order differential equation that exhibits limit cycle behavior, which makes it a valuable model for studying phenomena such as electrical circuits, cardiac rhythms, and neuronal firing patterns. The problem is governed by a parameter $\mu > 0$ and we consider its following formulation in the time interval $t \in [0, 3.5\mu]$

\begin{equation}
    \begin{cases}
        \dot{u}_1 = u_2,\\
        \dot{u}_2 = \mu (1 - u_1^2)u_2 - u_1,\\
        \mbf u(0) = (2, 0).\\
    \end{cases}
    \label{eq:vdp}
\end{equation}

This test case aims to assess the capabilities of our method on an autonomous periodic system. Indeed, the Van der Pol oscillator has two periodic solutions: a non-trivial periodic solution and an unstable zero solution. The parameter $\mu > 0$ determines the significance of the nonlinear component of the equation. When $\mu$ is large the system becomes stiff. We set $\Gamma = [10^2, 10^4]$ and consider the training dataset composed by a logarithmic discretization $\Gamma$ with $N_\mu=51$ points. The validation dataset is the interval midpoints of discretized $\Gamma$. The reference solutions are computed using the Radau method, with absolute and relative tolerances set to $10^{-12}$.

The data for the supervised training is smoothed with a Savitzky-Golay filter of window seven and order two. Data is randomly subsampled to have about one thousand points for each period of the timeseries. Fine-tuning is done by increasing the size of the unroll length starting from 20 up to 80 (exact values are a hyperparameter to be tuned) and time-stepping $\Delta t = 1 / 40$. The starting learning rate is reduced by a couple of orders of magnitude with respect to the supervised step.

Data is normalized using the following functions: parameter normalization: $\hat \mu = \log_{10} \mu$; state normalization: $\hat{\mbf u} = (u_1 / 2, u_2 / \mu)$; dynamics normalization: $\hat{\mbf f} = (f_1 / 5, f_2 / 10) / \mu$; time normalization: $\hat T = T / \mu$. The model is tested for $N_\mu^\text{test}=4$ values of $\mu \in [10^{2.01}, 10^{3.99}]$ spaced logarithmically (these values are chosen so that they are not present in the training nor the validation dataset). Moreover, we are also testing for time beyond what is present in the training and validation datasets, namely, we are integrating until a final time $\hat T = 5$. Figure~\ref{fig:rom-vdp} shows the results of applying the change of dynamics to this problem. On the left, the plot of the problem with the new dynamics shows smaller gradients. On the right, after the change of variable, we recover with good accuracy the original dynamics of the Van der Pol oscillator. In Table~\ref{tab:stiff-vdp-cost} we show a comparison of the computational cost. The number of right-hand side evaluations is one order of magnitude smaller for the explicit solver of the neural ODE w.r.t.\ the implicit solver of the FOM, this proves that the change of variables has significantly reduced the stiffness of the problem. Another interesting feature of our approach is that the ROM error is almost independent from $\mu$, on the other hand, the implicit solver exhibits a correlation among the two, suggesting that for larger values of $\mu$ computational gains may be larger. In this case, the $L^2$ error is computed only with respect to the first component of $\mbf u$.  Indeed, since $u_2$ features sharp peaks, small variations in their position make the integral error large and not representative of the accuracy of the solution. For this reason, we have introduced the $d_\text{peak}$ metric, which tries to measure the accuracy of the ROM in this regard. We stress that for this test case, the ROM is performing well also for unknown times, maintaining an accurate prediction of the period of the system. From the metrics in Table~\ref{tab:stiff-vdp-cost} it is also possible to notice that the main bottleneck in terms of accuracy for the ROM is the time mapping. Indeed, the model is very accurate in the $t_s$ parametrization.

In this case, we employ a fixed time-stepping in the online phase. This is a feasible approach only because the dynamics in the $t_s$ parameterization exhibit similar behavior across all considered values of $\mu$.

\begin{table}[t]
    \small
    \centering
    \begin{tabular}{llllllllll}
        $\mu$ & solver & tol & time [s] & \# fev & \# jev & \# lu & MSE${}_{t_s}$ & $L^2$ & $d_{\text{peak}}$ \\
        \hline
        \multirow{3}{*}{$10^{2.01}$} & ROM                    & $1/40$      & 0.033 & 1000 & 0 & 0 & 4.13e-3 & 0.03520 & 0.662\\
                                     & \multirow{2}{*}{Radau} & $10^{-2}$   & 0.046 & 2470 & 110 & 516 & -- & 0.01447 & 0.811\\
                                     &                        & $10^{-1.5}$ & 0.041 & 2342 & 105 & 482 & -- & 0.37252 & 14.59\\
        \hline
        \multirow{3}{*}{$10^{2.67}$} & ROM                    & $1/40$      & 0.033    & 1000 & 0 & 0 & 2.39e-3 & 0.02262 & 2.665\\
                                     & \multirow{2}{*}{Radau} & $10^{-2}$   & 0.056 & 3316 & 133 & 688 & -- & 0.00870 & 2.075\\
                                     &                        & $10^{-1.5}$ & 0.046 & 2591 & 118 & 576 & -- & 0.53100 & 181.5\\
        \hline
        \multirow{3}{*}{$10^{3.33}$} & ROM                    & $1/40$      & 0.033 & 1000 & 0 & 0 & 5.94e-3 & 0.03098 & 19.27\\
                                     & \multirow{2}{*}{Radau} & $10^{-2}$   & 0.072 & 3850 & 154 & 806 & -- & 0.00742 & 8.083\\
                                     &                        & $10^{-1.5}$ & 0.052 & 2992 & 137 & 684 & -- & 0.29987 & 655.2\\
        \hline
        \multirow{3}{*}{$10^{3.99}$} & ROM                    & $1/40$      & 0.033 & 1000 & 0 & 0 & 7.23e-3 & 0.04768 & 191.0\\
                                     & \multirow{2}{*}{Radau} & $10^{-2}$   & 0.078 & 4418 & 156 & 938 & -- & 0.04377 & 210.9\\
                                     &                        & $10^{-1.5}$ & 0.058 & 3337 & 148 & 768 & -- & 0.17724 & 2822\\
                                     
    \end{tabular}
    \caption{Test case 1: Van der Pol. Comparison of computational cost and accuracy for the Radau solver and the neural ODE based reduced order model (ROM) on the test dataset. Refer to Section~\ref{sec:numerical-results} for the definition of the metrics.}
    \label{tab:stiff-vdp-cost}
\end{table}


\subsection{Test case 2: OREGO problem}
\begin{figure}[!t]
    \centering
    \begin{minipage}{0.25\textwidth}
        \centering
        \includegraphics[trim={0 0.9cm 0 0.8cm},clip,height=3.0cm]{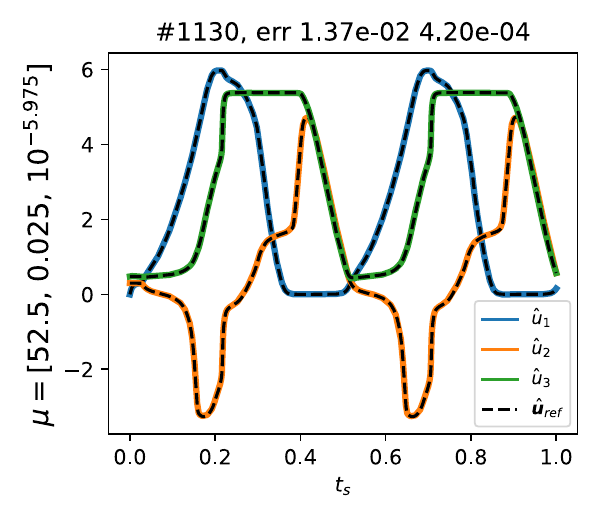}
    \end{minipage}%
    \begin{minipage}{0.25\textwidth}
        \centering
        \includegraphics[trim={1cm 0.9cm 0.4cm 0.8cm},clip,height=3.0cm]{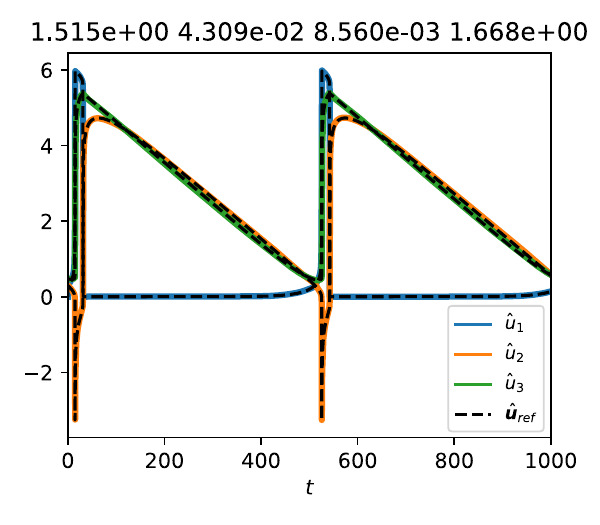}
    \end{minipage}\vrule%
    \begin{minipage}{0.25\textwidth}
        \centering
        \includegraphics[trim={0 0.9cm 0 0.8cm},clip,height=3.0cm]{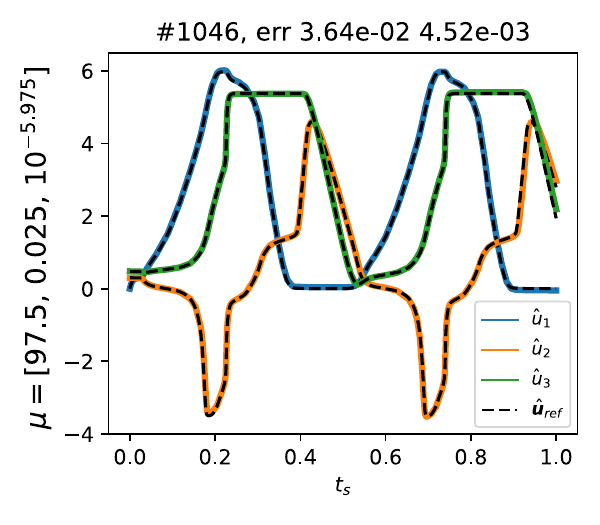}
    \end{minipage}%
    \begin{minipage}{0.25\textwidth}
        \centering
        \includegraphics[trim={1cm 0.9cm 0.4cm 0.8cm},clip,height=3.0cm]{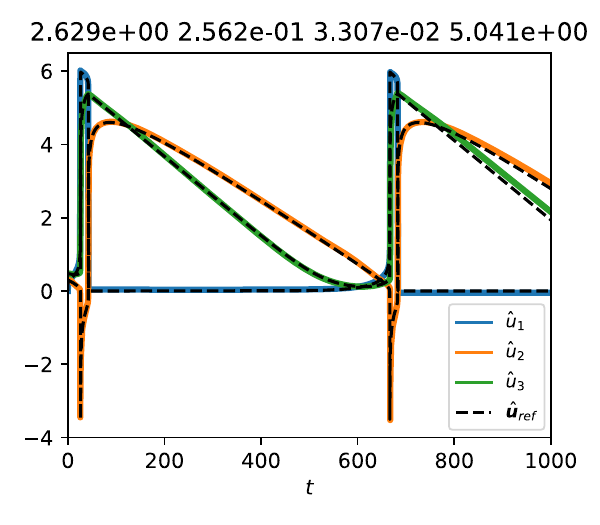}
    \end{minipage}

    \begin{minipage}{0.25\textwidth}
        \centering
        \includegraphics[trim={0 0.9cm 0 0.8cm},clip,height=3.0cm]{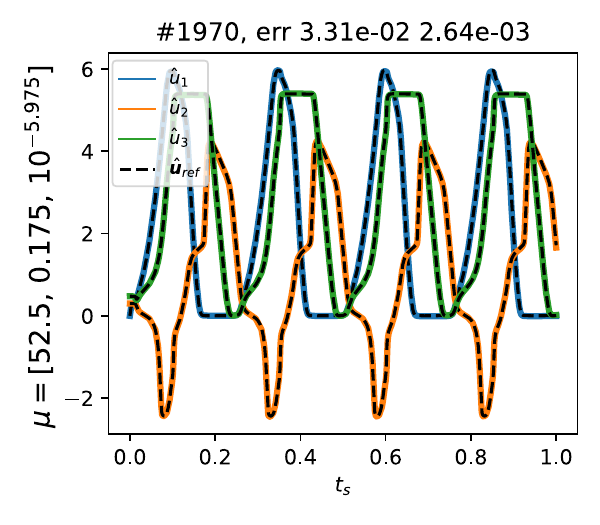}
    \end{minipage}%
    \begin{minipage}{0.25\textwidth}
        \centering
        \includegraphics[trim={1cm 0.9cm 0.4cm 0.8cm},clip,height=3.0cm]{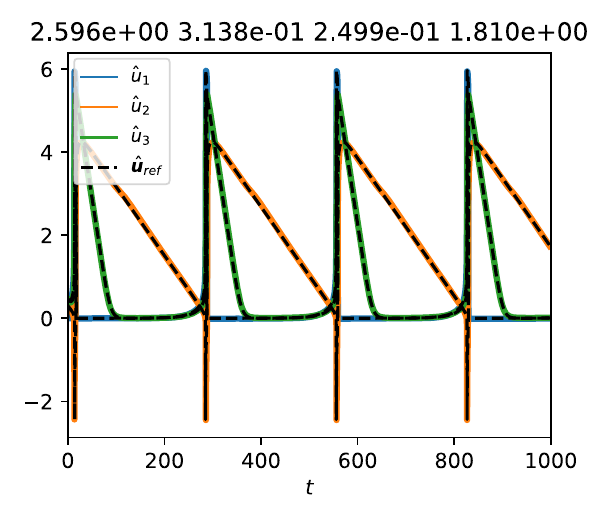}
    \end{minipage}\vrule%
    \begin{minipage}{0.25\textwidth}
        \centering
        \includegraphics[trim={0 0.9cm 0 0.8cm},clip,height=3.0cm]{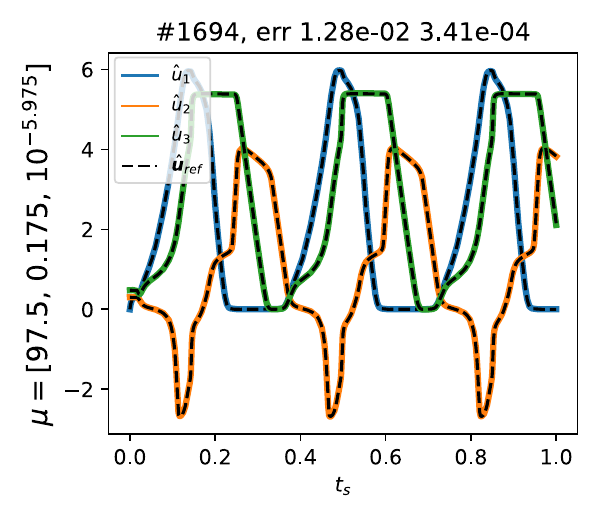}
    \end{minipage}%
    \begin{minipage}{0.25\textwidth}
        \centering
        \includegraphics[trim={1cm 0.9cm 0.4cm 0.8cm},clip,height=3.0cm]{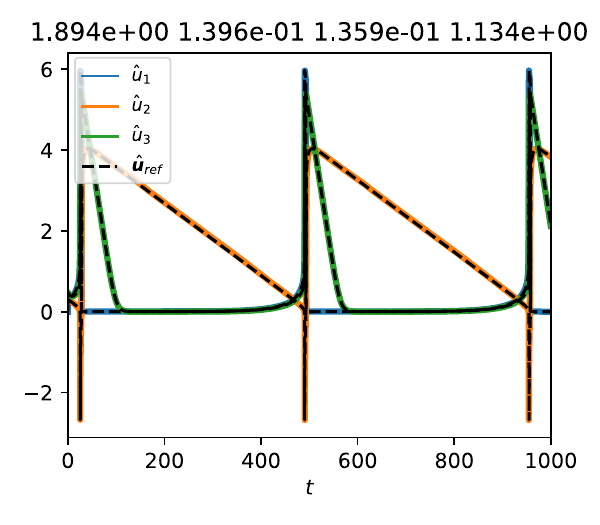}
    \end{minipage}

    \begin{minipage}{0.25\textwidth}
        \centering
        \includegraphics[trim={0 0.9cm 0 0.8cm},clip,height=3.0cm]{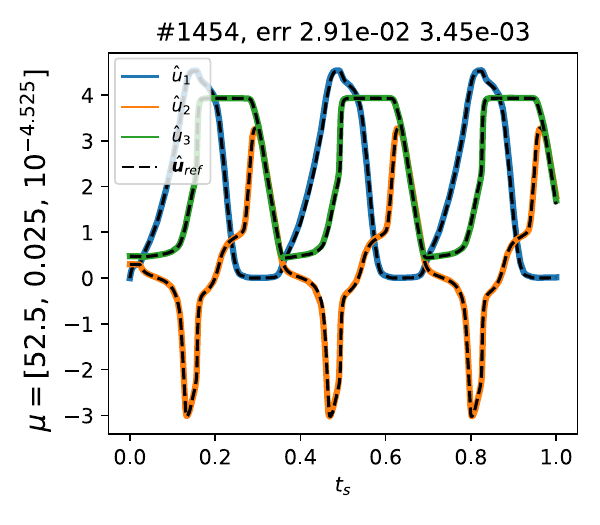}
    \end{minipage}%
    \begin{minipage}{0.25\textwidth}
        \centering
        \includegraphics[trim={1cm 0.9cm 0.4cm 0.8cm},clip,height=3.0cm]{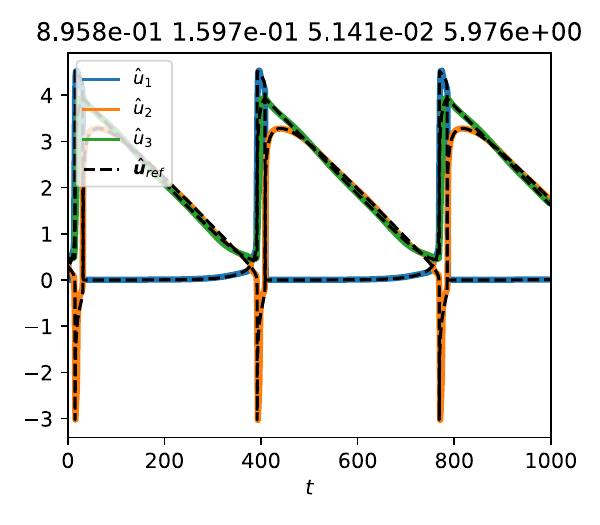}
    \end{minipage}\vrule%
    \begin{minipage}{0.25\textwidth}
        \centering
        \includegraphics[trim={0 0.9cm 0 0.8cm},clip,height=3.0cm]{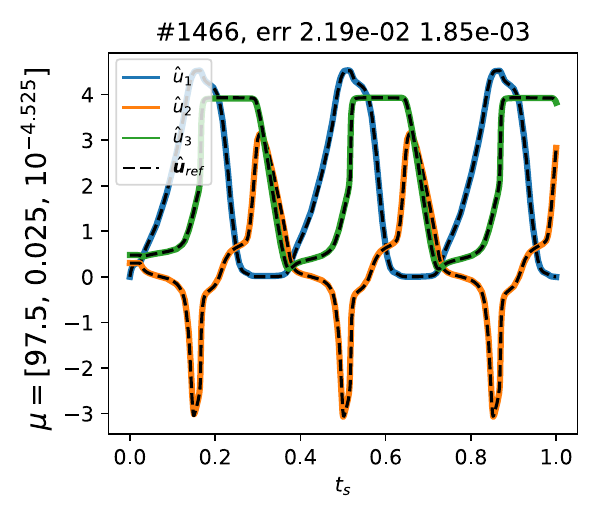}
    \end{minipage}%
    \begin{minipage}{0.25\textwidth}
        \centering
        \includegraphics[trim={1cm 0.9cm 0.4cm 0.8cm},clip,height=3.0cm]{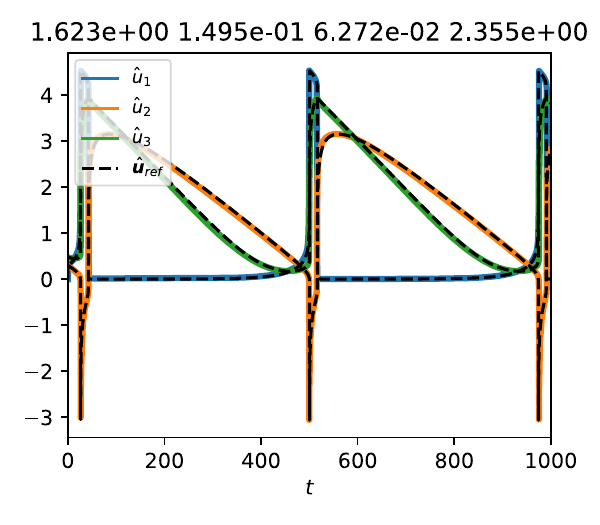}
    \end{minipage}

    \begin{minipage}{0.25\textwidth}
        \centering
        \includegraphics[trim={0 0 0 0.8cm},clip,height=3.35cm]{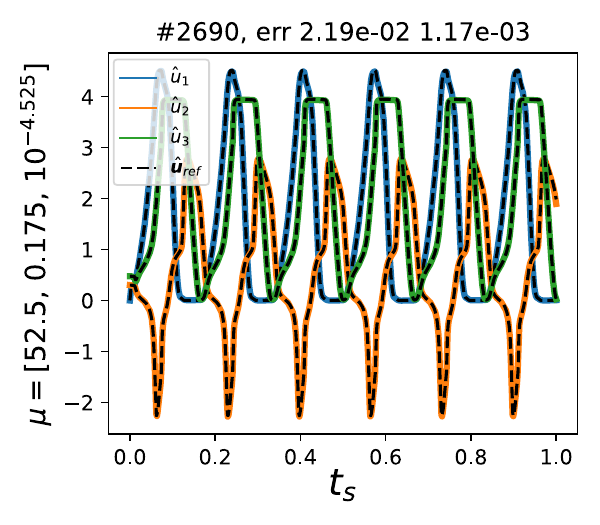}
    \end{minipage}%
    \begin{minipage}{0.25\textwidth}
        \centering
        \includegraphics[trim={1cm 0 0.4cm 0.8cm},clip,height=3.35cm]{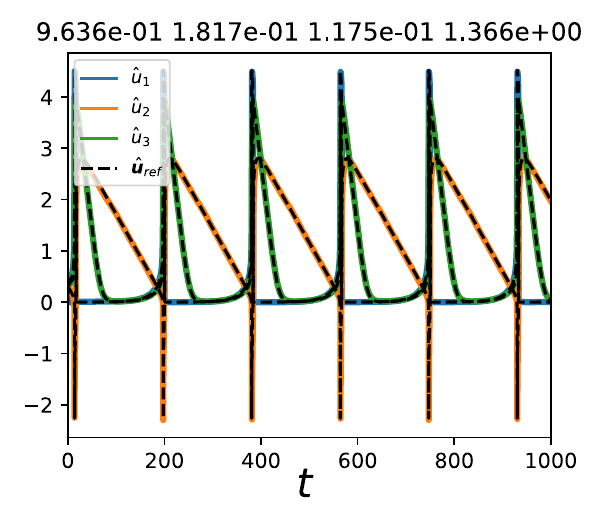}
    \end{minipage}\vrule%
    \begin{minipage}{0.25\textwidth}
        \centering
        \includegraphics[trim={0 0 0 0.8cm},clip,height=3.35cm]{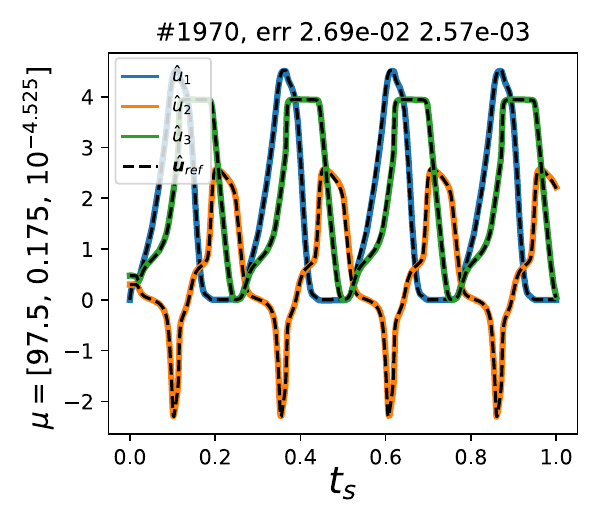}
    \end{minipage}%
    \begin{minipage}{0.25\textwidth}
        \centering
        \includegraphics[trim={1cm 0 0.4cm 0.8cm},clip,height=3.35cm]{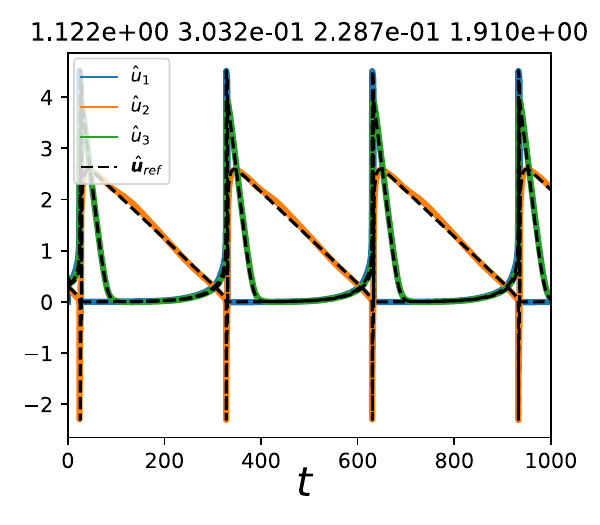}
    \end{minipage}

    \caption{Test case 2: OREGO. On the left, in colors, the neural ODE predictions in $t_s$. On the right, in colors, the prediction is mapped to the original time $t$. In black is the reference solution. The state scale is logarithmic.}
    \label{fig:rom-orego}
\end{figure}

The OREGO problem consists of the following stiff system of 3 non-linear ODEs defined in the time interval $t \in [0, 1000]$

\begin{equation}
    \begin{cases}
        \dot{u}_1 = \mu_1 (u_2 - u_1 u_2 + u_1 - \mu_3 u_1^2),\\
        \dot{u}_2 = \frac{1}{\mu_1} (-u_2 - u_1 u_2 + u_3), \\
        \dot{u}_3 = \mu_2 (u_1 - u_3), \\
        \mbf u(0) = (1, 2, 3).\\
    \end{cases}
    \label{eq:orego}
\end{equation}

Similarly to the previous test case, the problem is an autonomous system that features a periodic solution. However, the problem originates from the description of a chemical reaction, thus it is key for the concentration $\mbf u$ to stay positive. Moreover, the values of $\mbf u$ span six orders of magnitude. The considered parameter space is $\Gamma=[50, 100]\times[0.002, 0.02]\times[10^{-6}, 10^{-4}]$, which contains the value $\mbf \mu=(77.27, 0.161, 8.375 \cdot 10^{-5})$ usually employed in literature for this test case. The training dataset is built by subsampling a discretization of $\Gamma$ with a uniform grid of 6, 19, and 21 points in each direction, respectively. The validation dataset is based on a discretization of $\Gamma$ in the midpoints of the training discretization. Reference solutions are computed with a tolerance of $10^{-10}$. We apply the following normalizations: parameter normalization: $\hat{\mbf \mu} = (\mu_1 / 77.27, \mu_2 / 0.161, \log_{10} \mu_3 / 5)$; state normalization $\hat{\mbf u} = \log_{10} \mbf u$; dynamics normalization: $\hat{\mbf f} = \mbf f$. We test the performance of the ROM for eight parameters 
\begin{equation*}
    \mbf \mu^\text{test} \in \{52.5,97.5\} \times\{0.025,0.175\} \times\{10^{-5.975}, 10^{-4.525}\}.
\end{equation*}

\begin{figure}[!t] 
    \centering
    \includegraphics[width=0.8\textwidth]{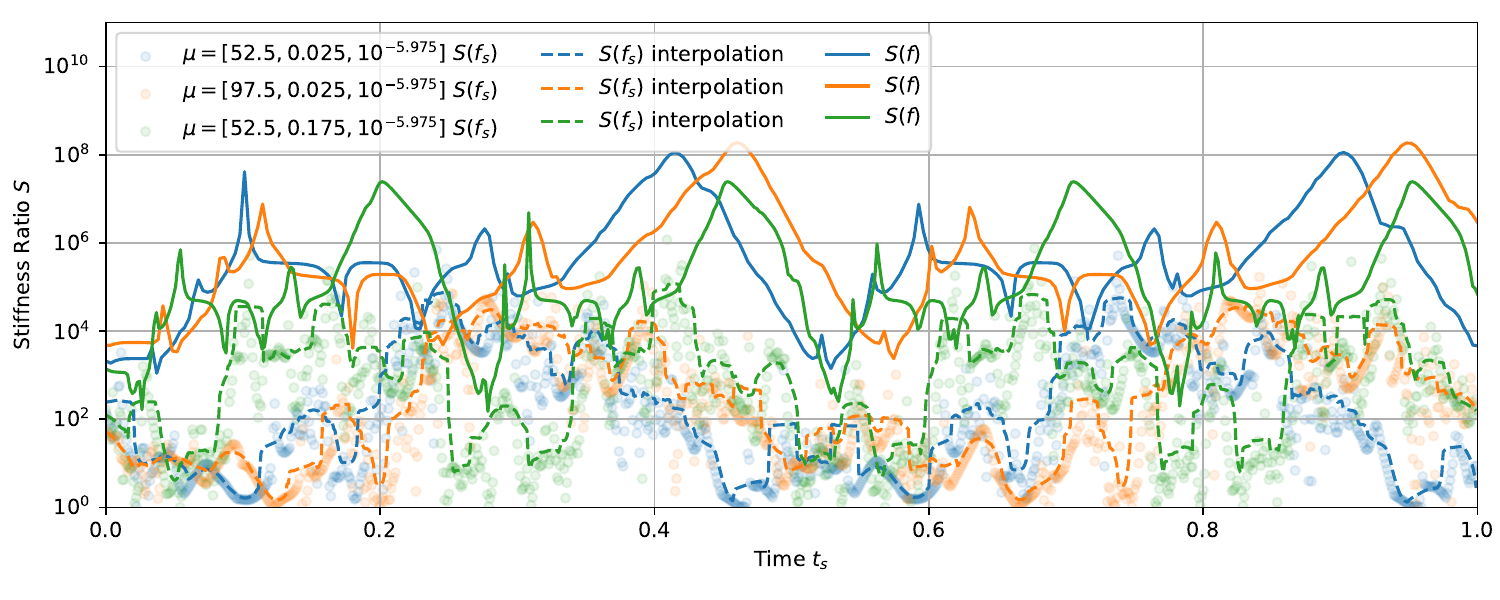}
    \caption{ Comparison of the stiffness ratio $S$ of OREGO system in the original and simplified dynamics.}
    \label{fig:orego-stiffness}
\end{figure}

Figure~\ref{fig:rom-orego} shows that the ROM is able to accurately follow the reference solution and track its period. A quantitative comparison with the Radau solver is reported in Table~\ref{tab:orego-cost}. We remark that a tolerance of $10^{-2}$ for the Radau solver did not produce a positive solution for all the problems in the test dataset. When this happens, the error is reported in bold in the table. The number of right-hand side evaluations of the explicit solver is less than half the one of the implicit one, showing that our method has made the system nonstiff. Also, the computational cost is consistently smaller than the one needed by the Radau solver, however, the latter has better accuracy. On the other hand, the distance of the peaks is comparable, proving that the ROM accurately captures the period of the system. In Figure~\ref{fig:orego-stiffness} we present a comparison of the stiffness ratio for the original system and the one with simplified dynamics. On average, the time reparametrization reduces the stiffness by a couple of orders of magnitude, showing the effectiveness of our technique.

\begin{table}[t]
    \small
    \centering
    \begin{tabular}{llllllllll}
        $\mbf \mu$ & solver & tol & time [s] & \# fev & \# jev & \# lu & MSE${}_{t_s}$ & $L^2$ & $d_{\text{peak}}$ \\
        \hline
        \multirow{3}{*}{$(52.5, 0.025, 10^{-5.975})$} & ROM                    & $2\cdot10^{-4}$ & 0.037 & 1130 & 0   & 0   & 1.20e-3 & 1.19e-2 & 1.668\\
                                                      & \multirow{2}{*}{Radau} & $10^{-3}$         & 0.065 & 3777 & 167 & 656 & --      & 1.31e-5 & 3.032\\
                                                      &                        & $10^{-1.5}$       & 0.059 & 3213 & 165 & 664  & -- & \textbf{0.4691} & 6.092\\
        \hline
        \multirow{2}{*}{$(97.5, 0.025, 10^{-5.975})$} & ROM   & $2\cdot10^{-4}$ & 0.034 & 1046 & 0   & 0   & 4.52e-3 & 7.75e-2 & 5.041\\
                                                      & \multirow{2}{*}{Radau} & $10^{-3}$         & 0.067 & 3445 & 159 & 600 & --      & 2.12e-5 & 4.632\\
                                                      &                        & $10^{-1.5}$       & 0.054 & 2815 & 141 & 556  & -- & \textbf{0.0099} & 6.574\\
        \hline
        \multirow{2}{*}{$(52.5, 0.175, 10^{-5.975})$} & ROM   & $2\cdot10^{-4}$ & 0.065 & 1970 & 0   & 0   & 2.64e-3 & 5.67e-2 & 1.810\\
                                                      & \multirow{2}{*}{Radau} & $10^{-3}$         & 0.111 & 6590 & 310 & 1100& --      & 8.47e-5 & 1.115\\
                                                      &                        & $10^{-1.5}$       & 0.095 & 5342 & 277 & 1068 & -- & 0.0079 & 1.432\\
        \hline
        \multirow{2}{*}{$(97.5, 0.175, 10^{-5.975})$} & ROM   & $2\cdot10^{-4}$ & 0.057 & 1694 & 0   & 0   & 1.41e-3 & 3.08e-2 & 1.134\\
                                                      & \multirow{2}{*}{Radau} & $10^{-3}$         & 0.074 & 4649 & 212 & 744 & --      & 1.43e-5 & 0.467\\
                                                      &                        & $10^{-1.5}$       & 0.068 & 3654 & 197 & 712  & -- & 0.0163 & 3.831\\
        \hline
        \multirow{2}{*}{$(52.5, 0.025, 10^{-4.525})$} & ROM   & $2\cdot10^{-4}$ & 0.048 & 1454 & 0   & 0   & 3.45e-3 & 5.14e-2 & 5.976\\
                                                      & \multirow{2}{*}{Radau} & $10^{-3}$         & 0.063 & 3895 & 183 & 666 & --      & 3.73e-5 & 2.885\\
                                                      &                        & $10^{-1.5}$       & 0.033 & 1814 & 92  & 364  & -- & \textbf{3.8838} & 464.3\\
        \hline
        \multirow{2}{*}{$(97.5, 0.025, 10^{-4.525})$} & ROM   & $2\cdot10^{-4}$ & 0.049 & 1466 & 0   & 0   & 1.85e-3 & 6.27e-2 & 2.355\\
                                                      & \multirow{2}{*}{Radau} & $10^{-3}$         & 0.062 & 3809 & 172 & 658 & --      & 8.83e-5 & 7.598\\
                                                      &                        & $10^{-1.5}$       & 0.037 & 2021 & 100 & 392  & -- & \textbf{2.6084} & 257.8\\
        \hline
        \multirow{2}{*}{$(97.5, 0.175, 10^{-4.525})$} & ROM   & $2\cdot10^{-4}$ & 0.088 & 2690 & 0   & 0   & 2.17e-3 & 1.17e-2 & 1.366\\
                                                      & \multirow{2}{*}{Radau} & $10^{-3}$         & 0.117 & 7148 & 319 & 1202& --      & 1.33e-5 & 1.100\\
                                                      &                        & $10^{-1.5}$       & 0.098 & 5433 & 292 & 1054 & -- & 0.0039 & 1.888\\
        \hline
        \multirow{2}{*}{$(52.5, 0.175, 10^{-4.525})$} & ROM   & $2\cdot10^{-4}$ & 0.059 & 1770 & 0   & 0   & 2.57e-3 & 2.28e-2 & 1.920\\
                                                      & \multirow{2}{*}{Radau} & $10^{-3}$         & 0.077 & 4781 & 217 & 814 & --      & 2.35e-5 & 1.102\\
                                                      &                        & $10^{-1.5}$       & 0.066 & 3635 & 196 & 712  & -- & 0.0101 & 0.992\\
             
    \end{tabular}
    \caption{Test case 2: OREGO. Comparison of computational cost and accuracy for the Radau solver and the neural ODE based reduced order model (ROM) on the test dataset. Refer to Section~\ref{sec:numerical-results} for the definition of the metrics. In bold, the error when the solution is not positive.}
    \label{tab:orego-cost}
\end{table}

\subsubsection{Parametric initial conditions}
\begin{figure}[!t]
    \centering
    \begin{minipage}{0.25\textwidth}
        \centering
        \includegraphics[trim={0 0 0 0},clip,height=3.0cm]{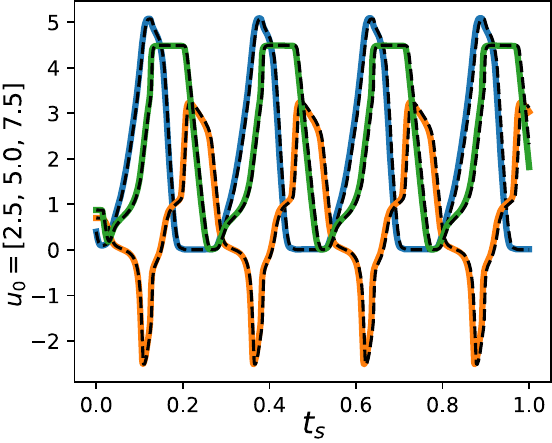}
    \end{minipage}%
    \begin{minipage}{0.25\textwidth}
        \centering
        \includegraphics[trim={0 0 0 0},clip,height=3.0cm]{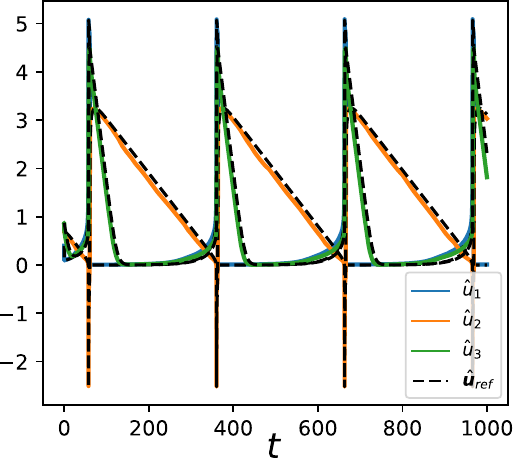}
    \end{minipage}\vrule%
    \begin{minipage}{0.25\textwidth}
        \centering
        \includegraphics[trim={0 0 0 0},clip,height=3.0cm]{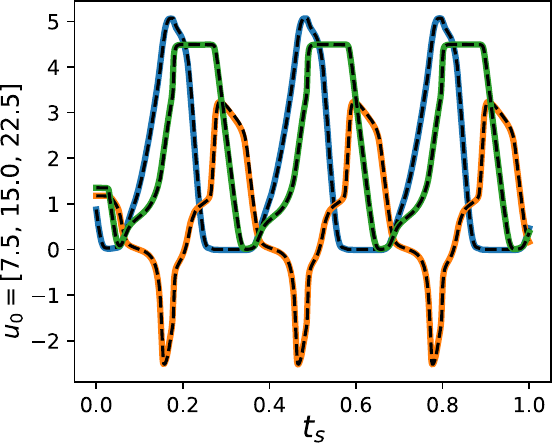}
    \end{minipage}%
    \begin{minipage}{0.25\textwidth}
        \centering
        \includegraphics[trim={0 0 0 0},clip,height=3.0cm]{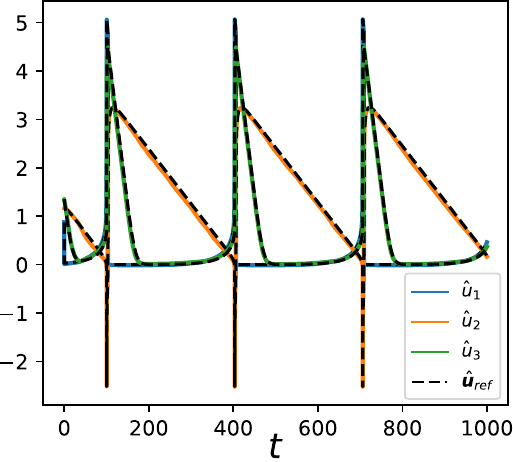}
    \end{minipage}
    \caption{Test case 2: OREGO, parametric initial conditions. On the left, in colors, the neural ODE predictions in $t_s$. On the right, in colors, the prediction is mapped to the original time $t$. In black is the reference solution. The state scale is logarithmic.}
    \label{fig:rom-orego-ic}
\end{figure}
To prove that the parameter $\mbf \mu$ is not limited to parametrizing the physical properties of the system, in this section we deal with the following variation of Eq.~\bref{eq:orego}
\begin{equation}
    \begin{cases}
        \dot{u}_1 = 77.27 (u_2 - u_1 u_2 + u_1 - 8.375\cdot10^{-6} u_1^2),\\
        \dot{u}_2 = 77.27^{-1} (-u_2 - u_1 u_2 + u_3), \\
        \dot{u}_3 = 0.161 (u_1 - u_3), \\
        \mbf u(0) = (\mu_1, \mu_2, \mu_3),\\
    \end{cases}
    \label{eq:orego-ic}
\end{equation}
where $\mbf \mu$ parameterizes the initial condition $\mbf u_0$. The training samples for the $\mbf \mu$ are in 
$$\Gamma = \{(1, 2, 3),(5, 10, 15), (10, 20, 30)\},$$
the normalization of the state space and the right-hand size is the same as before. The parameter $\mbf \mu$ is not normalized.
The problem is tested on $\mbf \mu^\text{test} \in \{(2.5, 5, 7.5), (7.5, 15, 22.5)\}$. The results are reported in Figure~\ref{fig:rom-orego-ic}. Our model is able to capture the period of the system and the different number of peaks present in the considered time window. The $L^2$-error is 0.0408 and 0.0256 in the first and second case, respectively.


\subsection{Test case 3: ROBER problem}

\begin{figure}[!t]
    \centering
    \begin{minipage}{0.5\textwidth}
        \centering
        \includegraphics[trim={0 0.9cm 0 0.8cm},clip,width=0.7\linewidth]{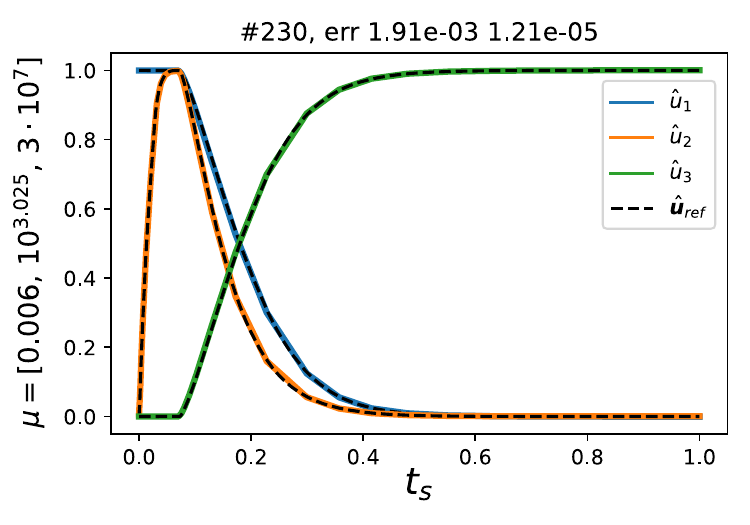}
    \end{minipage}%
    \begin{minipage}{0.5\textwidth}
        \centering
        \includegraphics[trim={0 0.9cm 0 0.8cm},clip,width=0.7\linewidth]{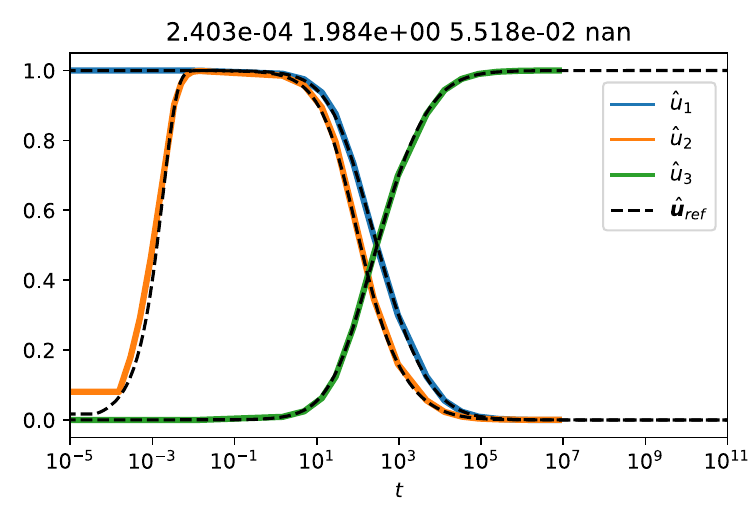}
    \end{minipage}
    \begin{minipage}{0.5\textwidth}
        \centering
        \includegraphics[trim={0 0.9cm 0 0.8cm},clip,width=0.7\linewidth]{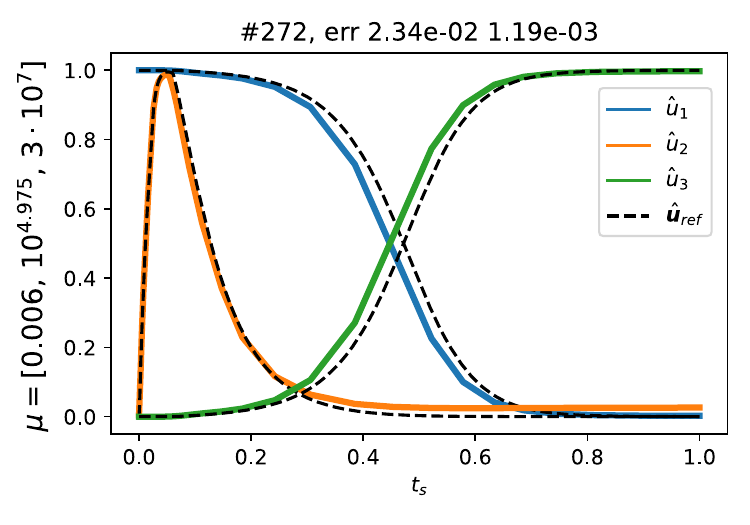}
    \end{minipage}%
    \begin{minipage}{0.5\textwidth}
        \centering
        \includegraphics[trim={0 0.9cm 0 0.8cm},clip,width=0.7\linewidth]{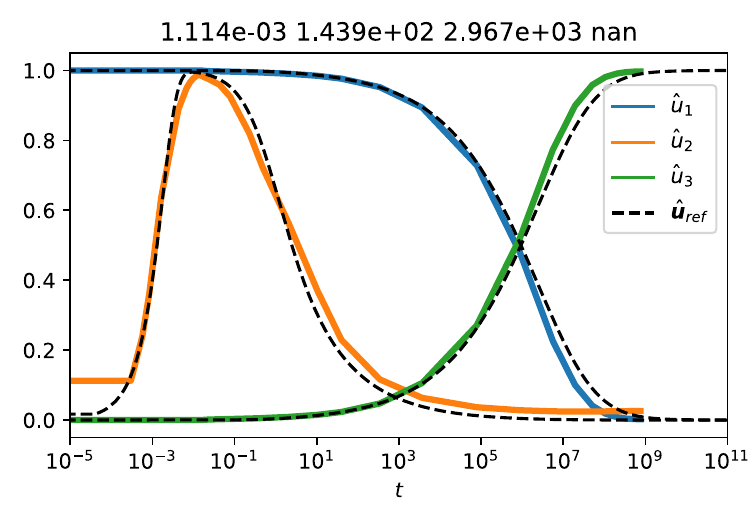}
    \end{minipage}
    \begin{minipage}{0.5\textwidth}
        \centering
        \includegraphics[trim={0 0.9cm 0 0.8cm},clip,width=0.7\linewidth]{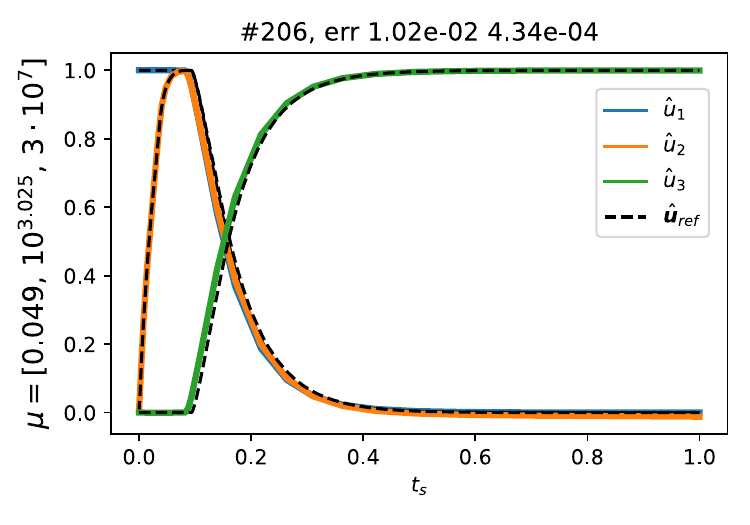}
    \end{minipage}%
    \begin{minipage}{0.5\textwidth}
        \centering
        \includegraphics[trim={0 0.9cm 0 0.8cm},clip,width=0.7\linewidth]{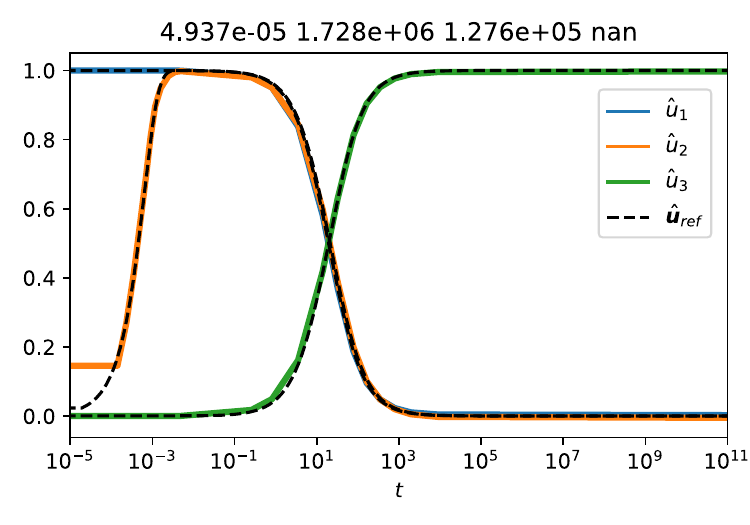}
    \end{minipage}
    \begin{minipage}{0.5\textwidth}
        \centering
        \includegraphics[trim={0 0 0 0.8cm},clip,width=0.7\linewidth]{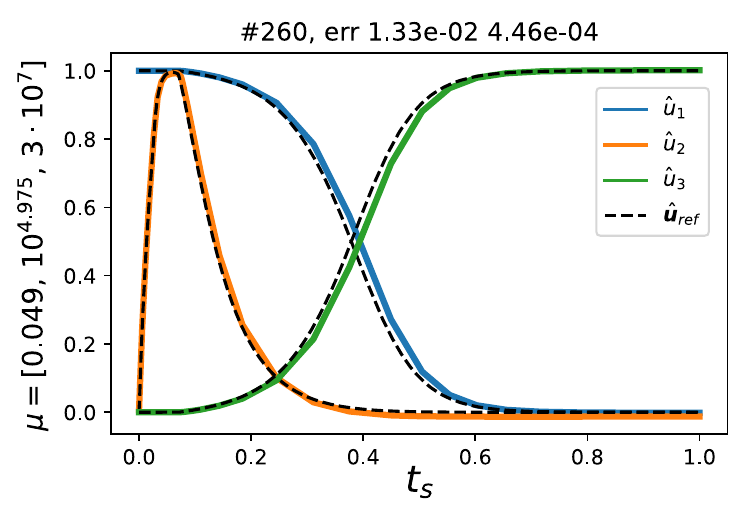}
    \end{minipage}%
    \begin{minipage}{0.5\textwidth}
        \centering
        \includegraphics[trim={0 0 0 0.8cm},clip,width=0.7\linewidth]{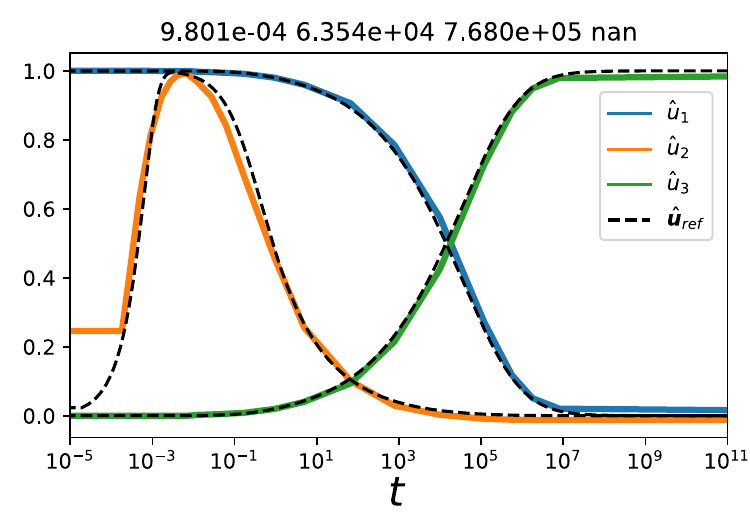}
    \end{minipage}
    \caption{Test case 3: ROBER. On the left, in colors, the neural ODE predictions in $t_s$. On the right, in colors, the prediction is mapped to the original time $t$. In black is the reference solution.}
    \label{fig:rom-rober}
\end{figure}

The ROBER problem is a stiff system of 3 non-linear ODEs defined as follows:

\begin{equation}
    \begin{cases}
        \dot{u}_1 = -\mu_1 u_1 + \mu_2 u_2 u_3,\\
        \dot{u}_2 = \mu_1 u_1 - \mu_2 u_2 u_3 - \mu_3 u_2^2, \\
        \dot{u}_3 = \mu_3 u_2^2,\\
        \mbf u(0) = (1, 0, 0).\\
    \end{cases}
    \label{eq:rober}
\end{equation}

The system describes the kinetics of an autocatalytic chemical reaction and is considered one of the most popular benchmarks for stiff solvers. When used in this context, the time integration window used is usually large, namely $t \in [0, 10^{11}]$ and the parameters are $\mbf \mu = [0.04, 10^4, 3\cdot10^7]$. Indeed, codes might fail if $u_2$ accidentally becomes negative, since it then tends to $-\infty$, causing overflow. This feature makes the ROBER problem ideal for testing the ability of our ROM to be accurate across different time scales. 

We set $\Gamma = [0.005, 0.05] \times [10^3, 10^5] \times \{3\cdot10^7\}$ and consider training dataset composed by a logarithmic discretization $\Gamma$ with 16 uniform points and 31 logarithmically spaced points in the first two directions, respectively. The validation dataset is the interval midpoints of the discretized $\Gamma$. The reference solutions are computed using the Radau method, with absolute and relative tolerances set to \texttt{rtol}$=10^{-10}$, \texttt{atol}$=10^{-14}$. 

Data is normalized using the following functions: $\hat{\mbf \mu} = (-\log_{10} \mu_1, \log_{10} \mu_2 / 4, \log_{10} \mu_3 / 7)$; state normalization $\hat{\mbf u} = (u_1, u_2 / 10^{(\hat{\mu_1} - \hat{\mu_3}) / 2}, u_3)$; dynamics normalization: $\hat{\mbf f} = \mbf f$. The model is tested for four values of the parameters 
\begin{equation*}
    \mbf \mu^\text{test} \in \{0.006, 0.49\} \times \{10^{3.025}, 10^{4.975}\}.
\end{equation*}

\begin{figure}[!t]
    \centering
    \begin{minipage}{0.5\textwidth}
        \centering
        \includegraphics[width=0.75\linewidth]{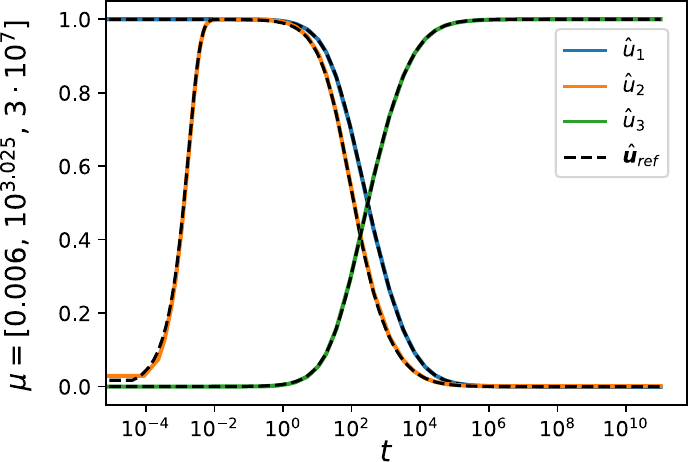}
    \end{minipage}%
    \begin{minipage}{0.5\textwidth}
        \centering
        \includegraphics[width=0.75\linewidth]{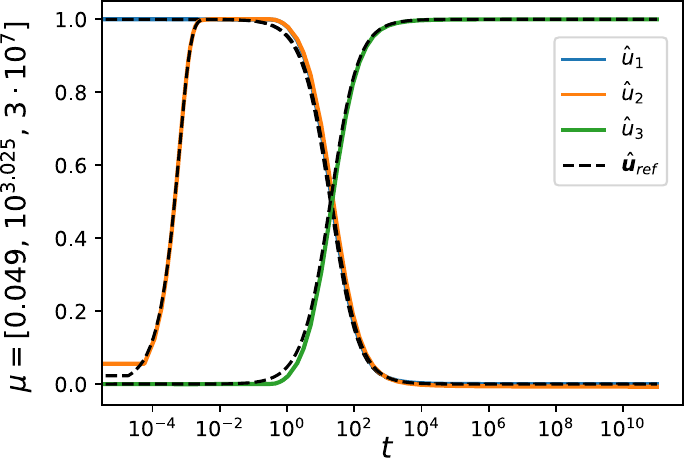}
    \end{minipage}
    \caption{Test case 3: ROBER. The prediction of the model with an integration tolerance of \texttt{tol}$=2\cdot 10^{-5}$ mapped to the original time $t$. The smaller tolerance results in more accurate predictions, see Figure~\ref{fig:rom-rober}.}
    \label{fig:rom-rober-lower-toll}
\end{figure}

Figure~\ref{fig:rom-rober} shows that the ROM is able to accurately follow the reference solution across the different time scales. A quantitative comparison with the Radau solver is reported in Table~\ref{tab:rober-cost}. We remark that a tolerance of \texttt{rtol}$=10^{-4}$, \texttt{atol}$=10^{-7}$ for the Radau solver is among the smallest pair of tolerances that produce a positive solution for all the problems in the test dataset. The number of right-hand side evaluations of the explicit solver is almost one order of magnitude smaller than the implicit one, showing that our method has made the system non-stiff. Moreover, the computational cost is consistently smaller than the one needed by the Radau solver. However, the latter achieves a far smaller error. We highlight that the ROM sometimes fails to reach the correct final time, namely, the map to $t$ sometimes stops at a time much smaller than $10^{11}$. Nevertheless, the ROM still captures the interesting part of the reaction, having issues only when the solution reaches a plateau. {Similarly, we observe that the model is not accurate for $t < 10^{-3}$ as it fails to follow the inflection point at $t\sim10^{-4}$ (cf.\ the right column of Figure~\ref{fig:rom-rober}). We show that these problems can be mitigated by reducing the tolerance of the RK45 solver to \texttt{tol}$=2\cdot 10^{-5}$, as shown in Figure~\ref{fig:rom-rober-lower-toll}.}

\begin{table}[t]
    \small
    \centering
    \begin{tabular}{lllllllll}
        $\mbf \mu$ & solver & rtol (atol) & time [s] & \# fev & \# jev & \# lu & MSE${}_{t_s}$ & MSE \\
        \hline
        \multirow{3}{*}{$(0.006, 10^{3.025}, 3 \cdot 10^7)$} & ROM   & $2\cdot10^{-4}$     & 0.008 &  230 & 0   & 0   & 1.21e-5 & 2.40e-4 \\
                                                             & \multirow{2}{*}{Radau} & $10^{-4} (10^{-7})$ & 0.019 & 927  & 44  & 210 & --      & 6.72e-6 \\
                                                             &                        & $10^{-3} (10^{-6})$ & 0.051 & 3055 & 132 & 798 & -- & 9.58e11 \\
        \hline
        \multirow{3}{*}{$(0.006, 10^{4.975}, 3 \cdot 10^7)$} & ROM   & $2\cdot10^{-4}$     & 0.009 &  272 & 0   & 0   & 1.19e-3 & 1.11e-3 \\
                                                             & \multirow{2}{*}{Radau} & $10^{-4} (10^{-7})$ & 0.021 & 1100 & 52  & 240 & --      & 5.74e-6 \\
                                                             &                        & $10^{-3} (10^{-6})$ & 0.016 & 790  & 53  & 194 & -- & 2.65e-5 \\
        \hline
        \multirow{3}{*}{$(0.049, 10^{3.025}, 3 \cdot 10^7)$} & ROM   & $2\cdot10^{-4}$     & 0.007 &  206 & 0   & 0   & 4.34e-4 & 4.97e-5 \\
                                                             & \multirow{2}{*}{Radau} & $10^{-4} (10^{-7})$ & 0.021 & 1081 & 46  & 256 & --      & 4.49e-6 \\
                                                             &                        & $10^{-3} (10^{-6})$ & 0.030 & 1787 & 73  & 406 & -- & 1.9626 \\
        \hline
        \multirow{3}{*}{$(0.049, 10^{4.975}, 3 \cdot 10^7)$} & ROM   & $2\cdot10^{-4}$     & 0.009 &  260 & 0   & 0   & 4.46e-4 & 9.80e-4 \\
                                                             & \multirow{2}{*}{Radau} & $10^{-4} (10^{-7})$ & 0.024 & 1121 & 55  & 248 & --      & 3.35e-6 \\
                                                             &                        & $10^{-3} (10^{-6})$ & 0.014 & 795  & 50  & 198 & -- & 1.50e-5 \\
    \end{tabular}
    \caption{Test case 3: ROBER. Comparison of computational cost and accuracy for the Radau solver and the neural ODE based reduced order model (ROM) on the test dataset. Refer to Section~\ref{sec:numerical-results} for the definition of the metrics.}
    \label{tab:rober-cost}
\end{table}

\subsection{Test case 4: E5 problem}

\begin{figure}[!t]
    \centering
    \begin{minipage}{0.5\textwidth}
        \centering
        \includegraphics[trim={0 0.9cm 0 0.77cm},clip,width=0.695\linewidth]{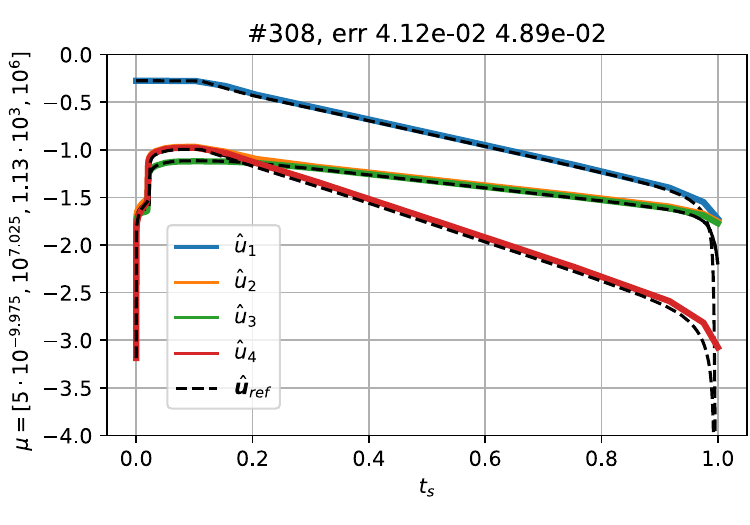}
    \end{minipage}%
    \begin{minipage}{0.5\textwidth}
        \centering
        \includegraphics[trim={0 0.9cm 0 0.77cm},clip,width=0.695\linewidth]{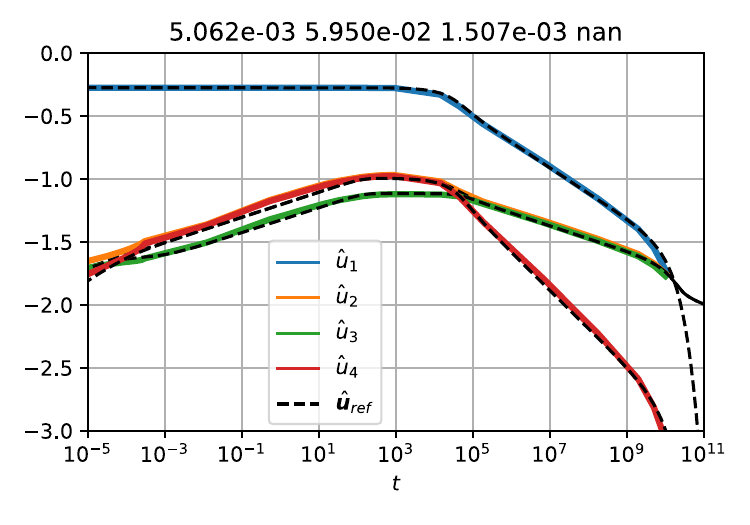}
    \end{minipage}
    \begin{minipage}{0.5\textwidth}
        \centering
        \includegraphics[trim={0 0.9cm 0 0.77cm},clip,width=0.695\linewidth]{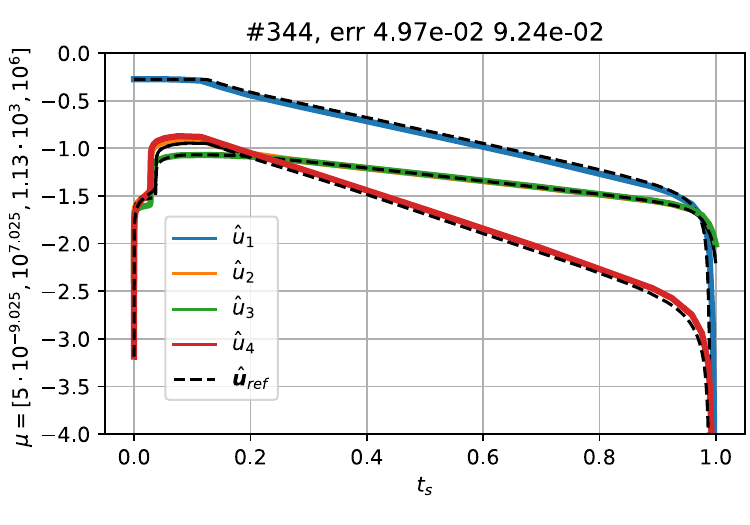}
    \end{minipage}%
    \begin{minipage}{0.5\textwidth}
        \centering
        \includegraphics[trim={0 0.9cm 0 0.77cm},clip,width=0.695\linewidth]{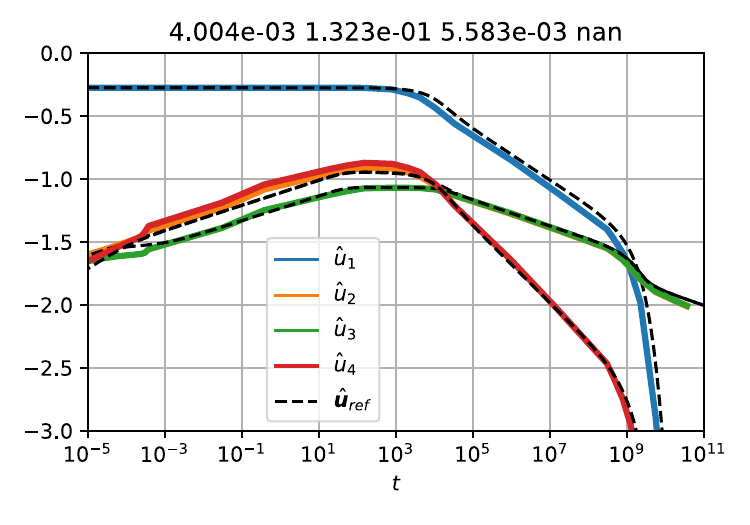}
    \end{minipage}
    \begin{minipage}{0.5\textwidth}
        \centering
        \includegraphics[trim={0 0.9cm 0 0.77cm},clip,width=0.695\linewidth]{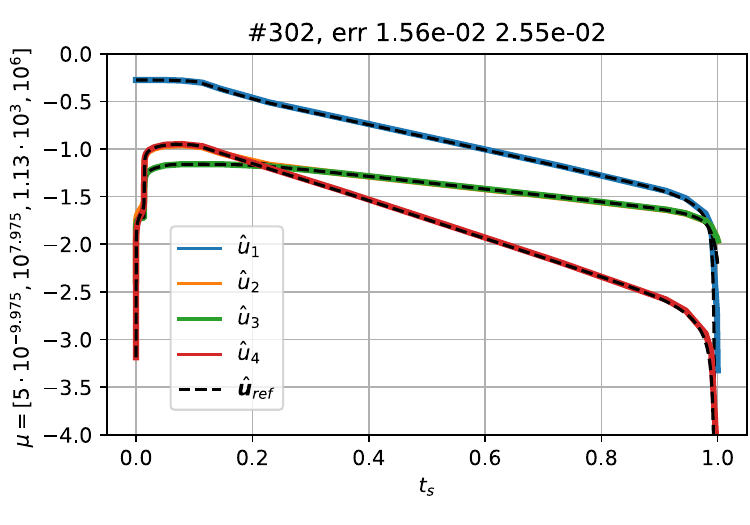}
    \end{minipage}%
    \begin{minipage}{0.5\textwidth}
        \centering
        \includegraphics[trim={0 0.9cm 0 0.77cm},clip,width=0.695\linewidth]{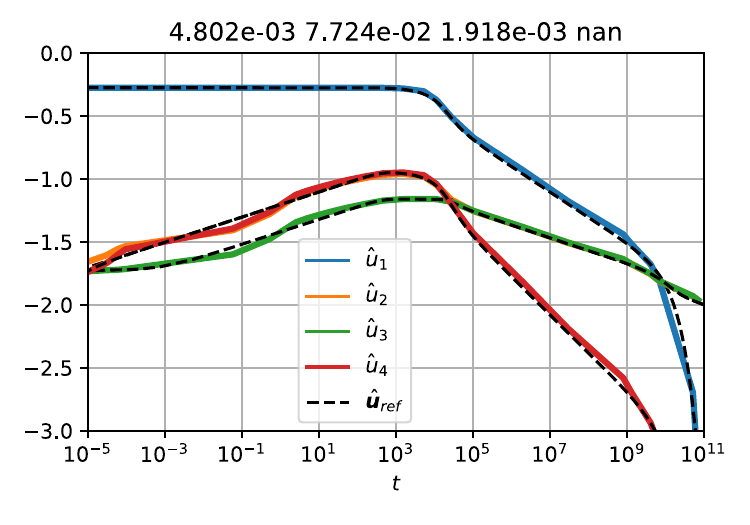}
    \end{minipage}
    \begin{minipage}{0.5\textwidth}
        \centering
        \includegraphics[trim={0 0 0 0.77cm},clip,width=0.695\linewidth]{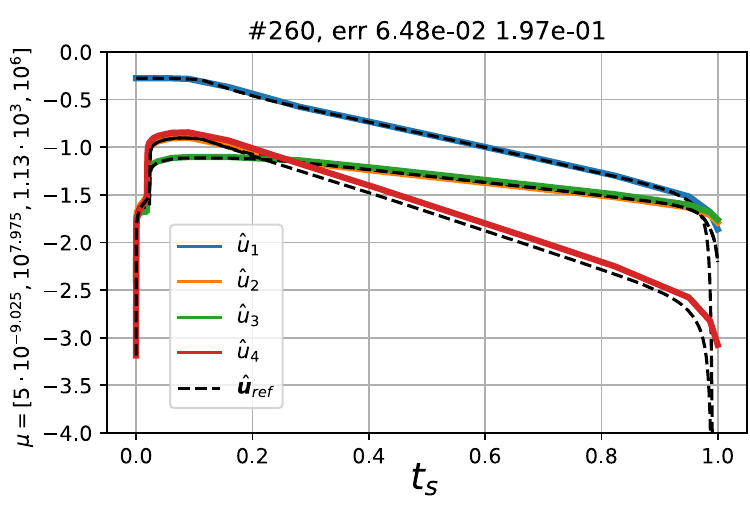}
    \end{minipage}%
    \begin{minipage}{0.5\textwidth}
        \centering
        \includegraphics[trim={0 0 0 0.77cm},clip,width=0.695\linewidth]{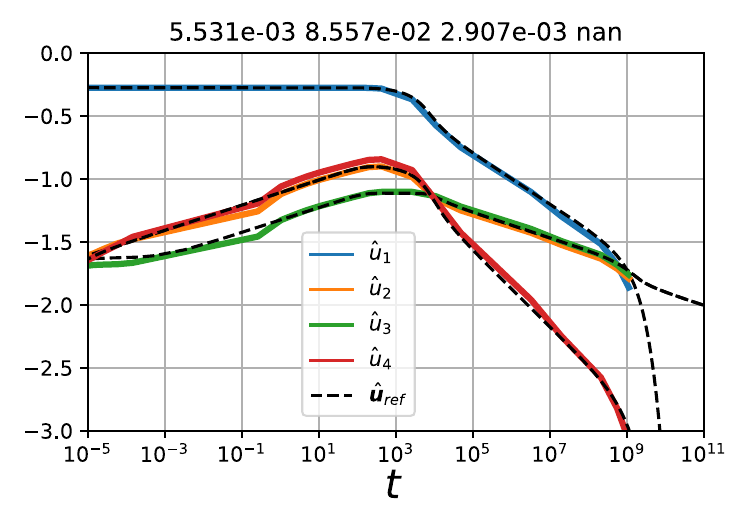}
    \end{minipage}
    \caption{Test case 4: E5. On the left, in colors, the neural ODE predictions in $t_s$. On the right, in colors, the prediction is mapped to the original time $t$. In black is the reference solution. The state scale is logarithmic.}
    \label{fig:rom-e5}
\end{figure}

The E5 problem consists of the following stiff system of 4 non-linear ordinary differential equations:

\begin{equation}
    \begin{cases}
        \dot{u}_1 = -\mu_1 u_1 + \mu_2 u_1 u_3,\\
        \dot{u}_2 = \mu_1 u_1 - \mu_3 \mu_4 u_2 u_3, \\
        \dot{u}_3 = \mu_1 u_1 - \mu_2 u_1 u_3 - \mu_3 \mu_4 u_2 u_3 + \mu_3 u_4,\\
        \dot{u}_4 = \mu_2 u_1 u_3 - \mu_3 u_4,\\
        \mbf u(0) = (1.76 \cdot 10^{-3}, 0, 0, 0).
    \end{cases}
    \label{eq:e5}
\end{equation}
The problem describes a chemical pyrolysis where $\mbf u$ represents the concentration of the reactants. The largely different rates of reaction that occur in the same system are the cause of stiffness. As a test problem, it is usually integrated in a large time window, namely, we choose the interval $t \in [0, 10^{11}]$. This means the problem features widely different scales in time and space. The problem is particularly challenging. Indeed, the formulation \bref{eq:e5} is affected by the cancellation of digits and it is recommended to use the relation $\dot{u}_3 = \dot{u}_2 - \dot{u}_4$ to compute $\mbf f(t, \mbf u)$. Moreover, it has been shown that a very small tolerance on the absolute scalar error must be used to obtain a reliable reference solution. Indeed, they are computed using the Radau method, with absolute and relative tolerances set to \texttt{atol}$=10^{-24}$, \texttt{rtol}$=10^{-10}$. 

We set $\Gamma = [5 \cdot 10^{-10}, 5 \cdot 10^{-9}] \times [10^7, 10^8] \times \{1.13\cdot10^3\} \times \{10^6\}$ 
which contains the parameter $\mbf \mu = (7.89\cdot 10^{-10}, 1.1\cdot 10^{7}, 1.13\cdot 10^{3}, 10^6)$, which is the usual choice for the benchmark.
The training dataset is built by a discretization $\Gamma$ with 11 logarithmically spaced points in the first two directions. The validation dataset is the interval midpoints of the discretized $\Gamma$. Data is normalized using the following functions: $\hat {\mbf \mu} = \log_{10} \mu $; state normalization $\hat{\mbf u} = \log_{10} \mbf u / 10$; dynamics normalization: 
\begin{equation*}
    \hat{\mbf f} = 
    \begin{cases}
        f_i & \text{ if } |f_i| < 2,\\
        \text{sgn}(f_i)(\log(|f_i| - 1) + 2) & \text{ otherwise. }
    \end{cases}
\end{equation*}
The insight behind this normalization is simple: the data features very long tails due to the presence of large gradients. Thus, we apply a logarithmic transformation to these large values. Suitable constants are added to connect with continuity with the linear part of the transformation.

The model is tested for four values of the parameters 
\begin{equation*}
    \mbf \mu^\text{test} \in \{5 \cdot 10^{-9.975}, 5 \cdot 10^{-9.025}\} \times \{10^{7.025}, 10^{7.975}\} \times \{1.13\cdot10^3\} \times \{10^6\}.
\end{equation*}
Figure~\ref{fig:rom-e5} shows that the ROM is able to accurately follow the reference solution across the different time and space scales. A quantitative comparison with the Radau solver is reported in Table~\ref{tab:e5-cost}. We remark that a tolerance of \texttt{rtol}$=10^{-6}$, \texttt{atol}$=10^{-20}$ for the Radau solver is among the smallest pair of tolerances that produce a positive and stable solution for all the problems in the test dataset. The number of right-hand side evaluations of the explicit solver is more than one order of magnitude smaller than the implicit one, showing that our method has indeed made the system nonstiff. The work precision tradeoff is particularly favorable for our ROM. Indeed, it has comparable accuracy at a much lower computational cost. The main limitation of our model is that it sometimes fails to reach the correct final time, namely, the map to $t$ sometimes stops at a time much smaller than $10^{11}$. However, the ROM still captures the interesting part of the reaction, having issues only when the concentration plummets to values close to zero ($<10^{-20}$).

\begin{table}[t]
    \small
    \centering
    \begin{tabular}{lllllllll}
        $\mbf \mu$ & solver & rtol (atol)  & time [s] & \# fev & \# jev & \# lu & MSE${}_{t_s}$ & MSE \\
        \hline
        \multirow{2}{*}{$(5\cdot10^{-9.975}, 10^{7.025}, $1.13e3$, 10^6)$} & ROM   & $2\cdot10^{-4}$      & 0.013 &  308 & 0   & 0   & 4.12e-2 & 5.06e-3 \\
                                                                                    & Radau & $10^{-6} (10^{-20})$ & 0.064 & 4119 & 93  & 416 & --      & 8.36e-4 \\
        \hline
        \multirow{2}{*}{$(5\cdot10^{-9.025}, 10^{7.025}, $1.13e3$, 10^6)$} & ROM   & $2\cdot10^{-4}$      & 0.015 &  344 & 0   & 0   & 4.97e-2 & 4.00e-3 \\
                                                                                    & Radau & $10^{-6} (10^{-20})$ & 0.065 & 4156 & 102 & 426 & --      & 1.84e-2 \\
        \hline
        \multirow{2}{*}{$(5\cdot10^{-9.975}, 10^{7.975}, $1.13e3$, 10^6)$} & ROM   & $2\cdot10^{-4}$      & 0.013 &  302 & 0   & 0   & 1.56e-2 & 4.80e-3 \\
                                                                                    & Radau & $10^{-6} (10^{-20})$ & 0.063 & 4003 & 103 & 430 & --      & 9.87e-4 \\
        \hline
        \multirow{2}{*}{$(5\cdot10^{-9.025}, 10^{7.975}, $1.13e3$, 10^6)$} & ROM   & $2\cdot10^{-4}$      & 0.011 &  260 & 0   & 0   & 6.48e-2 & 5.53e-3 \\
                                                                                    & Radau & $10^{-6} (10^{-20})$ & 0.068 & 4055 & 112 & 438 & --      & 1.88e-2 \\
    \end{tabular}
    \caption{Test case 4: E5. Comparison of computational cost and accuracy for the Radau solver and the neural ODE based reduced order model (ROM) on the test dataset. Refer to Section~\ref{sec:numerical-results} for the definition of the metrics.}
    \label{tab:e5-cost}
\end{table}


\subsection{Test case 5: POLLU problem}

\begin{figure}[!t]
    \centering
    $$\textnormal{Test case 5.1, } \mbf \mu = \mbf \mu^{\text{test} (1)}$$
    \begin{minipage}{0.5\textwidth}
        \centering 
        \includegraphics[width=0.95\linewidth]{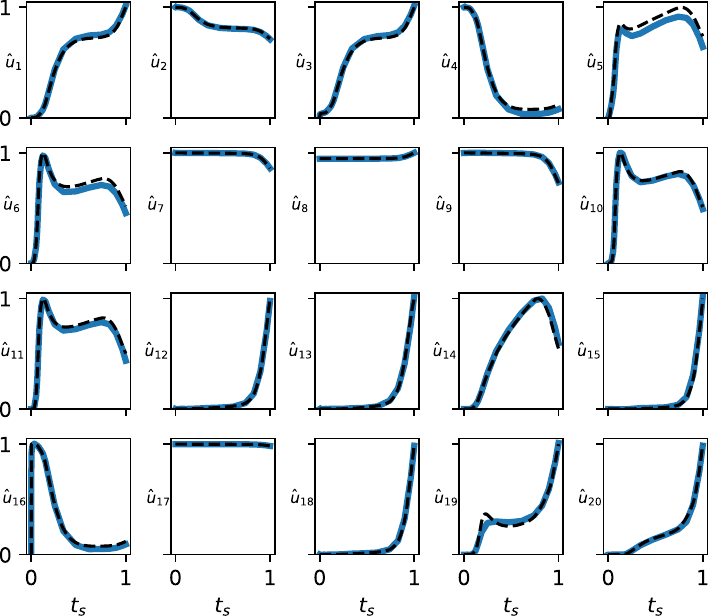}
    \end{minipage}\vrule
    \begin{minipage}{0.5\textwidth}
        \centering
        \includegraphics[width=0.95\linewidth]{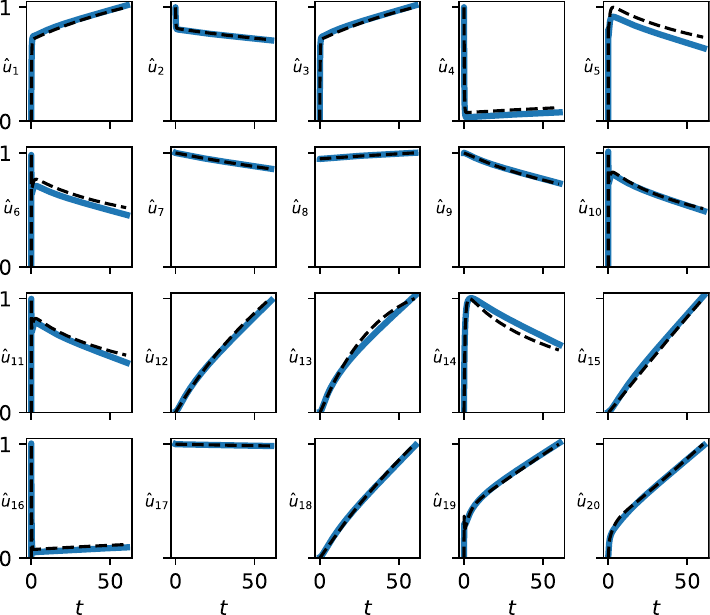}
    \end{minipage}
    $$\textnormal{Test case 5.2, } \mbf \mu = \mbf \mu^{\text{test} (2)}$$
    \begin{minipage}{0.5\textwidth}
        \centering
        \includegraphics[width=0.95\linewidth]{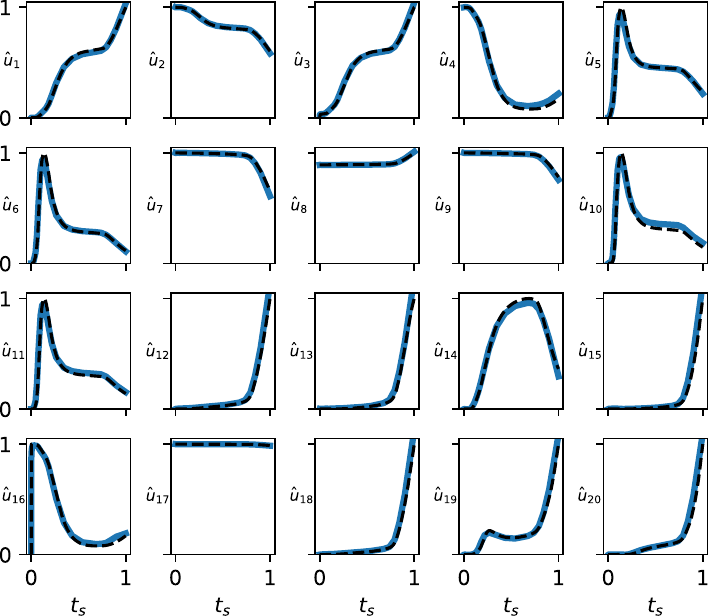}
    \end{minipage}\vrule
    \begin{minipage}{0.5\textwidth}
        \centering
        \includegraphics[width=0.95\linewidth]{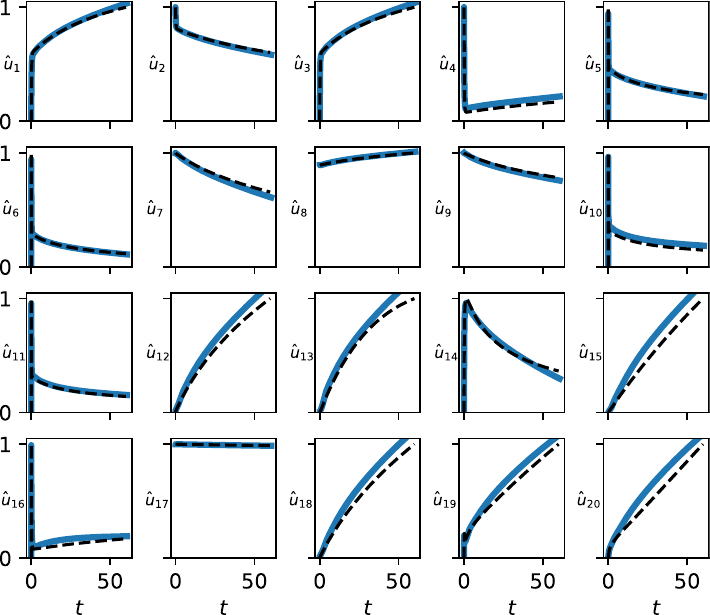}
    \end{minipage}

    \caption{Test case 5: POLLU. On the left, the neural ODE predictions in $t_s$. On the right, the prediction is mapped to the original time $t$. In blue (\includegraphics[width=0.4cm]{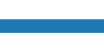}) the neural ODE prediction of $\hat u$, in black (\includegraphics[width=0.4cm]{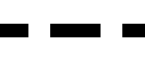}) the reference solution.}
    \label{fig:rom-pollu}
\end{figure}

To show the effectiveness of the proposed method on large state spaces, we test our methodology on the POLLU problem: as a stiff system of 20 non-linear ODEs. The system is the chemical reaction part of the air pollution model developed at the Dutch National Institute of Public Health and Environmental Protection. The problem features 25 reactions ($\mbf \mu \in \mathbb{R}^{25}$) and 20 reacting compounds ($\mbf u \in \mathbb{R}^{20}$). Namely, we have that the system is defined by

{\small
\begin{equation*}
    \mbf f(\mbf u; \mbf \mu) =
    \begin{bmatrix}
    -r_{1} - r_{10} - r_{14} - r_{23} - r_{24} + r_{2} + r_{3} + r_{9} + r_{11} + r_{12} + r_{22} + r_{25}\\
-r_{2} - r_{3} - r_{9} - r_{12} + r_{1} + r_{21}\\
-r_{15} + r_{1} + r_{17} + r_{19} + r_{22}\\
-r_{2} - r_{16} - r_{17} - r_{23} + r_{15}\\
-r_{3} + 2r_{4} + r_{6} + r_{7} + r_{13} + r_{20}\\
-r_{6} - r_{8} - r_{14} - r_{20} + r_{3} + 2r_{18}\\
-r_{4} - r_{5} - r_{6} + r_{13}\\
r_{4} + r_{5} + r_{6} + r_{7}\\
-r_{7} - r_{8}\\
-r_{12} + r_{7} + r_{9}\\
-r_{9} - r_{10} + r_{8} + r_{11}\\
r_{9}\\
-r_{11} + r_{10}\\
-r_{13} + r_{12}\\
r_{14}\\
-r_{18} - r_{19} + r_{16}\\
-r_{20}\\
r_{20}\\
-r_{21} - r_{22} - r_{24} + r_{23} + r_{25}\\
-r_{25} + r_{24}\\
    \end{bmatrix},
\end{equation*}}
where the reactions $\mbf r$ are connected to $\mbf u $ and $ \mbf \mu^{(0)}$ by rates reported in Table~\ref{tab:pollu-reaction-rates}. The initial condition is 
\begin{equation*}
    \mbf u_0 = (0, 0.2, 0, 0.04, 0, 0, 0.1, 0.3, 0.01, 0, 0, 0, 0, 0, 0, 0, 0.007, 0, 0, 0).
\end{equation*}

\begin{table}[!t]
    \small
    \centering
    \begin{tabular}{ccc|ccc|ccc}
        Reaction & $\mbf r$ & $\mbf \mu^{(0)}$ & Reaction & $\mbf r$ & $\mbf \mu^{(0)}$ & Reaction & $\mbf r$ & $\mbf \mu^{(0)}$ \\
        \hline
        1 & $\mu_{1}u_{1}      $  & 0.350e0 & 10 & $\mu_{10}u_{11}u_{1}$ & 0.900e4 & 19 & $\mu_{19}u_{16}     $ & 0.444e12 \\  
        2 & $\mu_{2}u_{2}u_{4} $  & 0.266e2 & 11 & $\mu_{11}u_{13}     $ & 0.220e-1 & 20 & $\mu_{20}u_{17}u_{6}$ & 0.124e4 \\      
        3 & $\mu_{3}u_{5}u_{2} $  & 0.120e5 & 12 & $\mu_{12}u_{10}u_{2}$ & 0.120e5 & 21 & $\mu_{21}u_{19}     $ & 0.210e1 \\ 
        4 & $\mu_{4}u_{7}      $  & 0.860e-3 & 13 & $\mu_{13}u_{14}     $ & 0.188e1 & 22 & $\mu_{22}u_{19}     $ & 0.578e1 \\
        5 & $\mu_{5}u_{7}      $  & 0.820e-3 & 14 & $\mu_{14}u_{1}u_{6} $ & 0.163e5 & 23 & $\mu_{23}u_{1}u_{4} $ & 0.474e-1 \\
        6 & $\mu_{6}u_{7}u_{6} $  & 0.150e5 & 15 & $\mu_{15}u_{3}      $ & 0.480e7 & 24 & $\mu_{24}u_{19}u_{1}$ & 0.178e4 \\
        7 & $\mu_{7}u_{9}      $  & 0.130e-3 & 16 & $\mu_{16}u_{4}      $ & 0.350e-3 & 25 & $\mu_{25}u_{20}     $ & 0.312e1 \\
        8 & $\mu_{8}u_{9}u_{6} $  & 0.240e5 & 17 & $\mu_{17}u_{4}      $ & 0.175e-1 &    & $                 $ &          \\
        9 & $\mu_{9}u_{11}u_{2}$  & 0.165e5 & 18 & $\mu_{18}u_{16}     $ & 0.100e9 &    & $                 $ &          \\
    \end{tabular}
    \caption{Reference reaction rates for the POLLU model that define $\mbf \mu^{(0)}$.}
    \label{tab:pollu-reaction-rates}
\end{table}

It is usually considered that the time interval $t \in [0,60]$ is representative of the behavior of the reactants sufficiently. The reference solution is obtained by using a tolerance of $10^{-10}$. Full details about the model can be found in \cite{verwer1994gauss}. The parameter space $\Gamma$ is built by changing the three components of $\mbf \mu$ that most impact the system. Namely, a preliminary sensibility analysis has shown that the components $4, 6$, and $14$ of $\mbf \mu$ are the most relevant, that is they produced the largest variation in the solution $\mbf u$ when perturbed. Hence, we define the parameter space as 
\begin{equation*}
    \Gamma = \bigtimes_{i \notin \{4, 6, 14\}} \{\mu_i^{(0)}\} \times \bigtimes_{i \in \{4, 6, 14\}} \left[\frac{1}{2}\mu_i^{(0)}, 2 \mu_i^{(0)}\right],
\end{equation*}
where $\mbf \mu^{(0)}$ is the reference parameter value used for the system, reported in Table~\ref{tab:pollu-reaction-rates}. The intervals are discretized with 16 uniform points in each coordinate direction. The validation dataset is made by the midpoints of the training datasets.

The model is tested for two parameters $\mbf \mu^{\text{test} (1)}, \mbf \mu^{\text{test} (2)}$ which are equal to $\mbf \mu^{(0)}$ apart from the three most relevant components, namely $\mu^{\text{test} (1)}_i = 0.525 \mu_i^{(0)}, \, \mu_i^{\text{test} (2)} = 1.975 \mu_i^{(0)}, \, i = 4, 6, 14$. Results are shown in Figure~\ref{fig:rom-pollu}. The ROM is able to follow the reference solution across all the reactions. A quantitative comparison with the Radau solver is reported in Table~\ref{tab:pollu-cost}. The number of right-hand side evaluations of the explicit solver is similar to the implicit one. On the other hand, each evaluation of the right-hand side, each evaluation of the Jacobian, and each LU decomposition is more costly than in the previous test cases since $N_u$ is larger. However, the cost of evaluating the neural network is similar to the previous test cases since the dimensions of the networks are comparable. Indeed, the cost of the ROM is still smaller than the cost of the Radau solver while keeping comparable accuracy. These results are of particular significance because they prove that as the state space gets larger our methodology becomes more cost effective.

\begin{table}[t]
    \small
    \centering
    \begin{tabular}{lllllllll}
        $\mbf \mu$ & solver & tol & time [s] & \# fev & \# jev & \# lu & MSE${}_{t_s}$ & MSE  \\
        \hline
        \multirow{3}{*}{$\mu_i = 0.525 \mu_i^{(0)}, \, i = 4, 6, 14$} & ROM   & $10^{-4}$ & 0.006 & 194 & 0  & 0  & 6.89e-4 & 9.80e-4 \\
                                                                  & \multirow{2}{*}{Radau} & $10^{-4}$ & 0.011 & 182 & 10 & 50 & --      & 3.81e-3\\
                                                                  &                        & $10^{-3}$ & 0.011 & 149 & 10 & 40 & --      & 0.2077\\
        \hline
        \multirow{3}{*}{$\mu_i = 1.975 \mu_i^{(0)}, \, i = 4, 6, 14$} & ROM   & $10^{-4}$ & 0.007 & 212 & 0  & 0  & 1.16e-3 & 1.69e-3 \\
                                                                  & \multirow{2}{*}{Radau} & $10^{-4}$ & 0.017 & 268 & 14 & 70 & --      & 9.03e-4 \\ 
                                                                  &                        & $10^{-3}$ & 0.014 & 203 & 12 & 54 & --      & 4.63e-3\\
    \end{tabular}
    \caption{Test case 5: POLLU. Comparison of computational cost and accuracy for the Radau solver and the neural ODE based reduced order model (ROM) on the test dataset. Refer to Section~\ref{sec:numerical-results} for the definition of the metrics.}
    \label{tab:pollu-cost}
\end{table}

\section{Conclusions}\label{sec:conclusions}

In this work, we have developed a novel methodology to tackle {some kinds of} stiffness in neural ODEs. The approach hinges on a suitable time reparametrization of the system that significantly reduces the stiffness, making it possible to efficiently apply explicit solvers to the neural ODEs.
The construction of the time reparametrization is general and completely data-driven, making it ideal for an application to ROMs.

By leveraging the intrinsic capabilities of neural ODEs to handle continuous-time data, we capture the dynamics of stiff systems effectively. The accuracy, robustness, and efficiency of the methodology were tested in five famous benchmark test cases from the literature. Namely, when applied to periodic systems, the model generalizes to data well beyond the time training interval. We attribute this characteristic primarily to the autonomous nature of the neural ODE and the time reparametrization model. When applied to chemical equations, our model always maintained the positivity of the concentration, even when applied to problems suffering from numerical instability. Moreover, the work/precision tradeoff of our models often proved comparable or favorable when confronted with a state-of-the-art Runge-Kutta implicit solver of the Radau II A kind. In a few cases where this was not observed, our model provided a computationally cheaper alternative for approximate solutions while preserving concentration positivity, which is not possible with the Radau solver for small tolerances.

In conclusion, our work highlights the promising potential of neural ODEs in creating efficient and accurate reduced-order models for stiff ODEs. {Stiffness characterized by reasonably large eigenvalues but very high frequency remains a challenge for our approach, and it is an active topic of research.} Future research is directed towards improving the accuracy of the time mapping that, as of right now, stands as the major accuracy bottleneck of the ROM. We also aim to test this framework to handle the semi-discrete formulation of partial differential equations and other higher-dimensional systems. Moreover, techniques like model distillation, network pruning, and better hyperparameter tuning could further reduce the computational cost of our approach.


\paragraph*{Declaration of competing interests.} 
The authors declare that they have no known competing financial interests or personal relationships that could have appeared to influence the work reported in this paper.

\paragraph*{Acknowledgments.}

MC has been partially funded by PRIN2020 n. 20204LN5N5 funded by the Italian Ministry of Universities and Research (MUR). MC is a member of INdAM-GNCS.

\paragraph*{CRediT authorship contribution statement.}
MC: Conceptualization, Methodology, Software, Validation, Formal analysis, Investigation, Data curation, Visualization, Writing - Original draft.
JSH: Conceptualization, Methodology, Resources, Writing - Review and editing, Supervision, Project administration.

\bibliographystyle{abbrv}
\bibliography{references}
\end{document}

%% file: tikz/fnn.tex
\colorlet{myred}{orange}
\colorlet{myblue}{blue!80!black}
\colorlet{mygreen}{green!60!black}
\colorlet{myorange}{orange!70!red!60!black}
\colorlet{mydarkred}{red!30!black}
\colorlet{mydarkblue}{blue!40!black}
\colorlet{mydarkgreen}{green!30!black}

\tikzset{
  >=latex, 
  node/.style={thick,circle,draw=myblue,minimum size=22,inner sep=0.5,outer sep=0.6},
  node in/.style={node,green!20!black,draw=myblue!30!black,fill=myblue!25},
  node hidden/.style={node,blue!20!black,draw=black,fill=black!10},
  node convol/.style={node,orange!20!black,draw=myorange!30!black,fill=myorange!20},
  node out/.style={node,red!20!black,draw=myred!30!black,fill=myred!20},
  connect/.style={thick,black}, 
  connect arrow/.style={-{Latex[length=4,width=3.5]},thick,black,shorten <=0.5,shorten >=1},
  node 1/.style={node in}, 
  node 2/.style={node hidden},
  node 3/.style={node out}
}
\def\nstyle{int(\lay<\Nnodlen?min(2,\lay):3)} 

\begin{tikzpicture}[x=2.2cm,y=1.4cm]
  \message{^^JNeural network without text}
  \readlist\Nnod{4,5,5,5,3} 
  
  \message{^^J  Layer}
  \foreachitem \N \in \Nnod{ 
    \def\lay{\Ncnt} 
    \pgfmathsetmacro\prev{int(\Ncnt-1)} 
    \message{\lay,}
    \foreach \i [evaluate={\y=\N/2-\i; \x=\lay; \n=\nstyle;}] in {1,...,\N}{ 
      
      \node[node \n] (N\lay-\i) at (\x,\y) {};
      
      \ifnum\lay>1 
        \foreach \j in {1,...,\Nnod[\prev]}{ 
          \draw[connect,white,line width=1.2] (N\prev-\j) -- (N\lay-\i);
          \draw[connect] (N\prev-\j) -- (N\lay-\i);
        }
      \fi 
      
    }
  }

  \node[node,red!20!black,draw=myred!30!black,fill=myred!20] (c) at (1,0) {};
  \node[node,red!20!black,draw=myred!30!black,fill=myred!20] (c) at (1,-1) {};
  \node[node,red!20!black,draw=myred!30!black,fill=myred!20] (c) at (1,-2) {};

  \draw[->,myred,line width=0.5mm] (5.2, 0.5) -> (5.8, 0.5);
  \draw[->,myred,line width=0.5mm] (5.2,-0.5) -> (5.8,-0.5);
  \draw[->,myred,line width=0.5mm] (5.2,-1.5) -> (5.8,-1.5);
  \node[fill=black!20,text=orange,rotate=270] at (5.45,-0.5) {\Large normalization};
  \node[node,red!20!black,draw=myred!30!black,fill=myred!20] (c) at (6, 0.5) {};
  \node[node,red!20!black,draw=myred!30!black,fill=myred!20] (c) at (6,-0.5) {};
  \node[node,red!20!black,draw=myred!30!black,fill=myred!20] (c) at (6,-1.5) {};
  \node[text=orange] at (6.65,-0.5) {\Huge $\tdev{\hat{\mbf u}}{t}$};
  
  \node[text=orange] at (-0.5,-1) {\LARGE $\hat{\mbf u}(0)$};
  \draw[->,myred,line width=0.5mm] (-0.1,-1) -> (0.75,-1);

  \node[text=orange] at (-0.5,-4) {\LARGE $\hat{\mbf u}(T)$};
  \draw[->,myred,line width=0.5mm] (0.75,-4) -> (-0.1,-4);

  \draw[->,myblue,line width=0.5mm] (-0.25,1) -> (0.75,1);
  \draw[dashed,->,myred,line width=0.5mm] (6.65,-1.2) -- (6.65,-4) -- (1,-4) -- (1,-2.5);
  \node[fill=black!20] at (5,-4) {\Large ODE Solver $\displaystyle \int \mathrm{d}t$ };

  \node[myblue] at (-0.5,1) {\Huge $\mbf \mu$};
  \draw[->,myblue,line width=0.5mm] (-0.25,1) -> (0.75,1);
  \node[fill=black!20,text=blue,rotate=270] at (0.25,1) {\large normalization};

\end{tikzpicture}

%% file: tikz/reparam_workflow.tex
\tikzset{every picture/.style={line width=0.75pt}} 

\begin{tikzpicture}[x=0.75pt,y=0.75pt,yscale=-1,xscale=1]

\draw  [draw opacity=0][fill={rgb, 255:red, 155; green, 155; blue, 155 }  ,fill opacity=1 ] (72.34,60.56) -- (72.34,82.64) -- (76.35,82.64) -- (68.34,97.36) -- (60.32,82.64) -- (64.33,82.64) -- (64.33,60.56) -- cycle ;
\draw (73, 220) node  {\includegraphics[width=110pt]{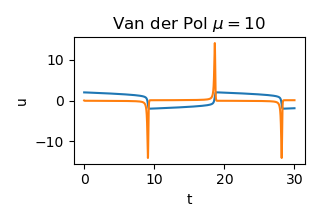}};
\draw (310, 117) node  {\includegraphics[width=110pt]{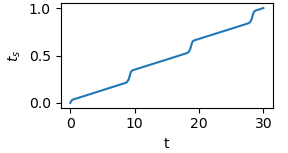}};
\draw (310, 233) node  {\includegraphics[width=110pt]{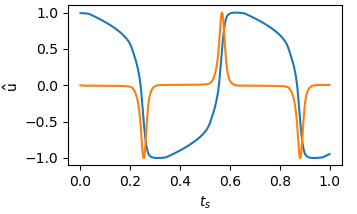}};
\draw (553, 239) node  {\includegraphics[width=110pt]{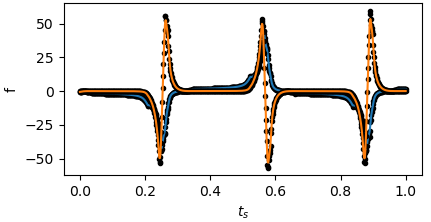}};
\draw (430, 70) node  {\includegraphics[width=50pt]{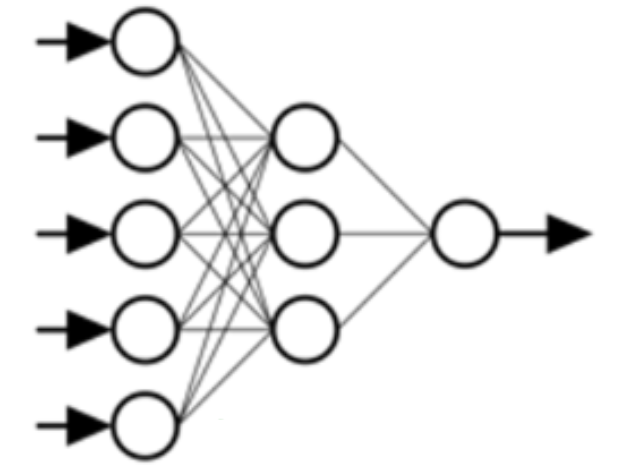}};
\draw (490, 145) node  {\includegraphics[width=50pt]{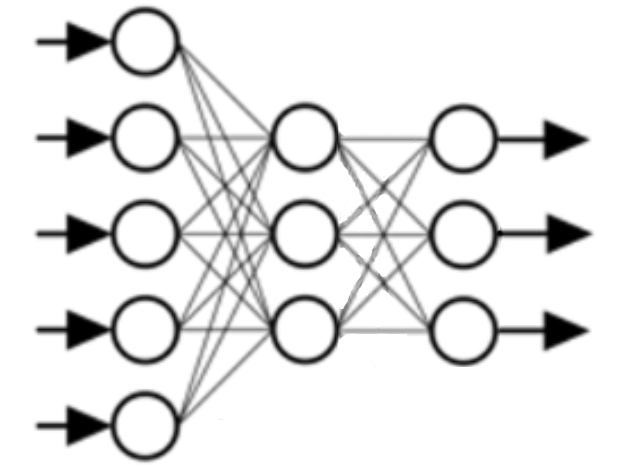}};
\draw   (0.00, 0.00) -- (150.00, 0.00) -- (150.00, 55.00) -- (00.00, 55.00) -- cycle ;
\draw   (0.00, 100.00) -- (150.00, 100.00) -- (150.00, 270.00) -- (0.00, 270.00) -- cycle ;
\draw  [draw opacity=0][fill={rgb, 255:red, 155; green, 155; blue, 155 }  ,fill opacity=1 ] (160.29,184.42) -- (181.06,184.42) -- (181.06,180.17) -- (194.91,188.68) -- (181.06,197.2) -- (181.06,192.94) -- (160.29,192.94) -- cycle ;
\draw   (235.0,42.57) -- (388.0,42.57) -- (388.0,277.02) -- (235.0,277.02) -- cycle ;

\draw  [draw opacity=0][fill={rgb, 255:red, 155; green, 155; blue, 155 }  ,fill opacity=1 ] (408.55,104.22) -- (429.32,104.22) -- (429.32,99.96) -- (443.17,108.48) -- (429.32,116.99) -- (429.32,112.74) -- (408.55,112.74) -- cycle ;
\draw  [draw opacity=0][fill={rgb, 255:red, 155; green, 155; blue, 155 }  ,fill opacity=1 ] (527.34,169.23) -- (527.34,147.15) -- (523.34,147.15) -- (531.35,132.43) -- (539.36,147.15) -- (535.36,147.15) -- (535.36,169.23) -- cycle ;
\draw   (470.0,15) -- (645.0,15) -- (645.0,119) -- (470.0,119) -- cycle ;
\draw   (470.0,173) -- (645.0,173) -- (645.0,276.92) -- (470.0,276.92) -- cycle ;
\draw  [fill={rgb, 255:red, 155; green, 155; blue, 155}  ,fill opacity=0.5 ] (0.00,0.00) -- (150.00, 0.00) -- (150.00, 25.00) -- (00.00, 25.00) ;
\draw  [fill={rgb, 255:red, 74; green, 144; blue, 226 }  ,fill opacity=0.5 ] (0.00, 100.00) -- (150.00, 100.00) -- (150.00, 125.00) -- (0.00, 125.00) -- cycle ;
\draw  [fill={rgb, 255:red, 74; green, 144; blue, 226 }  ,fill opacity=0.5 ] (235.0,42.57) -- (388.0,42.57) -- (388.0,67.57) -- (235.0,67.57) -- cycle ;
\draw  [fill={rgb, 255:red, 74; green, 144; blue, 226 }  ,fill opacity=0.5 ] (235.0,159.6) -- (388.0,159.6) -- (388.0,184.6) -- (235.0,184.6) -- cycle ;
\draw  [fill={rgb, 255:red, 245; green, 166; blue, 35 }  ,fill opacity=0.5 ] (470.0,15) -- (645.0,15) -- (645.0,40.0) -- (470.0,40.0) -- cycle ;
\draw  [fill={rgb, 255:red, 74; green, 144; blue, 226 }  ,fill opacity=0.5 ] (470.0, 173.0) -- (645.0,173.0) -- (645.0, 198.0) -- (470.0, 198.0) -- cycle ;

\draw  [draw opacity=0][fill={rgb, 255:red, 155; green, 155; blue, 155 }  ,fill opacity=1 ] (408.55,224.22) -- (429.32,224.22) -- (429.32,219.96) -- (443.17,228.48) -- (429.32,236.99) -- (429.32,232.74) -- (408.55,232.74) -- cycle ;

\draw (45, 5.8) node [anchor=north west][inner sep=0.75pt]   [align=left] {Stiff ODE};
\draw (4.74,32.45) node [anchor=north west][inner sep=0.75pt]   [align=left] {$\displaystyle \dot{\mbf u}\ =\ \mbf f( t,\ \mbf u( t) ;\ \mbf \mu )$};
\draw (81.3,70.74) node [anchor=north west][inner sep=0.75pt]  [font=\footnotesize,color={rgb, 255:red, 74; green, 74; blue, 74 }  ,opacity=1 ] [align=left] {Implicit solver};
\draw (17, 105) node [anchor=north west][inner sep=0.75pt]   [align=left] {Discretized solution};
\draw (5,127) node [anchor=north west][inner sep=0.75pt]   [align=left] {State: $\displaystyle \mbf u\left( t^{0}\right) ,\ ...,\ \mbf u\left( t^{n}\right)$};
\draw (5,147) node [anchor=north west][inner sep=0.75pt]   [align=left] {Time: $\displaystyle t^{0} ,\ ...,\ t^{n}$};
\draw (153.00, 125.00) node [anchor=north west][inner sep=0.75pt]  [font=\footnotesize,color={rgb, 255:red, 74; green, 74; blue, 74 }  ,opacity=1 ] [align=left] {Induced change \\of variable $\displaystyle t_{s}$};
\draw (153.00, 150.00) node [anchor=north west][inner sep=0.75pt]  [color={rgb, 255:red, 74; green, 74; blue, 74 }  ,opacity=1 ] [align=left] {$\displaystyle t_{s}\left( t^{i}\right) \ =\ \frac{i}{n}$};
\draw (163.00, 205.00) node [anchor=north west][inner sep=0.75pt]  [font=\footnotesize,color={rgb, 255:red, 74; green, 74; blue, 74 }  ,opacity=1 ] [align=left] {Normalize \\ state $\hat u$};
\draw (240,47) node [anchor=north west][inner sep=0.75pt]   [align=left] {Time reparametrization};
\draw (270,165) node [anchor=north west][inner sep=0.75pt]   [align=left] {New dynamics};
\draw (395,193.86) node [anchor=north west][inner sep=0.75pt]  [font=\footnotesize,color={rgb, 255:red, 74; green, 74; blue, 74 }  ,opacity=1 ] [align=left] {Interpolate};
\draw (390, 5) node [anchor=north west][inner sep=0.75pt]  [font=\footnotesize,color={rgb, 255:red, 74; green, 74; blue, 74 }  ,opacity=1 ] [align=left] {Normalize\\ and learn\\ $\displaystyle \mathcal{NN}_t \approx\dot{t}(\hat u)$};

\draw (450, 50) node [anchor=north west][inner sep=0.75pt]  [font=\small,color={rgb, 255:red, 74; green, 74; blue, 74 }  ,opacity=1 ] [align=left] {$\displaystyle \dot{t}$};
\draw (389, 55) node [anchor=north west][inner sep=0.75pt]  [font=\small,color={rgb, 255:red, 74; green, 74; blue, 74 }  ,opacity=1 ] [align=left] {$\displaystyle \hat{u}$};
\draw (389, 75) node [anchor=north west][inner sep=0.75pt]  [font=\small,color={rgb, 255:red, 74; green, 74; blue, 74 }  ,opacity=1 ] [align=left] {$\displaystyle \hat{\mu}$};

\draw (545.58,133) node [anchor=north west][inner sep=0.75pt]  [font=\footnotesize,color={rgb, 255:red, 74; green, 74; blue, 74 }  ,opacity=1 ] [align=left] {Normalize and learn\\ $\displaystyle \mathcal{NN} \approx \mbf f_{s}( t_{s} , \mbf u ; \mbf \mu )$};
\draw (505, 120) node [anchor=north west][inner sep=0.75pt]  [font=\small,color={rgb, 255:red, 74; green, 74; blue, 74 }  ,opacity=1 ] [align=left] {$\displaystyle f_{s}$};
\draw (450, 130) node [anchor=north west][inner sep=0.75pt]  [font=\small,color={rgb, 255:red, 74; green, 74; blue, 74 }  ,opacity=1 ] [align=left] {$\displaystyle \hat{\mbf u}$};
\draw (450, 150) node [anchor=north west][inner sep=0.75pt]  [font=\small,color={rgb, 255:red, 74; green, 74; blue, 74 }  ,opacity=1 ] [align=left] {$\displaystyle \hat{\mbf \mu}$};

\draw (475, 50) node [anchor=north west][inner sep=0.75pt]   [align=left] {\small
$\displaystyle \! \dot{\mbf u_s}(t_s; \mbf \mu)  = \mathcal{NN}( t_{s}, \mbf u_s( t_{s}) ; \mbf \mu )$\\[2mm]
\small
$\displaystyle \! t(t_s; \mbf \mu) = \! \int_{\mbf u_s(t_s^{0})}^{\mbf u_s(t_s)} \!\!\!\!\!\! \mathcal{NN}_t(\mbf u_s; \mbf \mu) \text{d}\mbf u_s$
};
\draw (500.00, 21) node [anchor=north west][inner sep=0.75pt]   [align=left] {Non-stiff surrogate};
\draw (520, 179) node [anchor=north west][inner sep=0.75pt]   [align=left] {Estimated $\displaystyle f_{s}$};

\end{tikzpicture}